\documentclass{amsart}
\usepackage{amssymb,amsthm,amsmath}
\usepackage[all]{xy}
\usepackage[latin1]{inputenc}
\usepackage{graphics}
\input{diagxy}
\usepackage{tikz}
\usetikzlibrary{shapes}
\usetikzlibrary{decorations.markings}

\newcommand{\myemph}[1]{\textbf{#1}}
\newcommand{\intervalsin}{\textbf{Int}}

\newcommand{\tensor}{\otimes}
\newcommand{\cotensor}{\pitchfork}
\newcommand{\cocomp}{\star}
\newcommand{\cocat}{\textnormal{\textbf{Cocat}}}
\newcommand{\CC}{\mathcal{C}}
\newcommand{\DD}{\mathcal{D}}
\newcommand{\EE}{\mathcal{E}}
\newcommand{\KK}{\mathcal{K}}
\newcommand{\sets}{\textnormal{\textbf{Set}}}
\newcommand{\topcat}{\textnormal{\textbf{Top}}}
\newcommand{\groupoids}{\textnormal{\textbf{Gpd}}}
\newcommand{\cat}{\textnormal{\textbf{Cat}}}
\newcommand{\chplus}{\textnormal{\textbf{Ch}}_{0\leq}}
\newcommand{\rmod}[1]{#1\textnormal{\textbf{-mod}}}

\newcommand{\iso}{\cong}
\newcommand{\bottom}{\bot}
\newtheorem{theorem}{Theorem}[section]
\newtheorem{lemma}[theorem]{Lemma}
\newtheorem{proposition}[theorem]{Proposition}
\newtheorem{corollary}[theorem]{Corollary}

\def\twocell{\to/=>/}
\tikzset{->-/.style={decoration={
  markings,
  mark=at position .5 with {\arrow{>}}},postaction={decorate}}}

\theoremstyle{definition}
\newtheorem{definition}[theorem]{Definition}
\newtheorem{example}[theorem]{Example}

\theoremstyle{remark}
\newtheorem*{remark}{Remark}
\begin{document}
%% 
%% TITLE INFORMATION
%% 
\author[M. A. Warren]{Michael A. Warren}
\email{mwarren@math.ias.edu}
\date{\today}
\address{School of Mathematics\\Institute for Advanced Study\\Einstein
  Dr., Princeton, New Jersey\\USA}
\title[representable intervals]{a characterization of
  representable intervals}
%%
%% MAKETITLE, TOC, ETC
%%
\begin{abstract}
  In this note we provide a characterization, in terms of additional
  algebraic structure, of those strict intervals
  (certain cocategory objects) in a symmetric monoidal closed category $\EE$ that
  are representable in the sense of inducing on $\EE$ the structure of
  a finitely bicomplete 2-category.  Several examples and
  connections with the homotopy theory of 2-categories are also discussed.
\end{abstract}
\maketitle
%%
%% Introduction
%%
\section*{Introduction}

Approached from an abstract perspective, a
fundamental feature of the category of spaces which enables the
development of homotopy theory is the presence of an object $I$
with which the notions of \emph{path} and \emph{deformation} thereof are
defined.  When dealing with topological spaces, $I$ is most
naturally taken to be the closed unit interval $[0,1]$, but there are
other instances where the homotopy theory of a category is determined in an
appropriate way by an interval object $I$.  For example,
the simplicial interval $I=\Delta[1]$ determines --- in a sense
clarified by the recent work of Cisinski
\cite{Cisinski:PCMTH} --- the classical model structure
on the category of simplicial sets and the infinite dimensional sphere
$J$ is correspondingly related to the quasi-category model structure
studied by Joyal \cite{Joyal:TQCA}.  Similarly, the
category $\mathbf{2}$ gives rise to the natural model
structure --- in which the weak equivalences are categorical
equivalences, the fibrations are isofibrations and the cofibrations
are functors injective on objects --- on the category $\cat$ of small categories
\cite{Joyal:SSCS}.  This model structure is, moreover, well-behaved
with respect to the usual 2-category structure on $\cat$ (it is a
\emph{model $\cat$-category} in the sense of \cite{Lack:HA2M}).  One
special property of the category $\mathbf{2}$, which is in part
responsible for these facts, is that it is a cocategory in $\cat$.  

In this paper we study, with a view
towards homotopy theory, one (abstract) notion of \emph{strict interval object} ---
namely, a cocategory with object of coobjects the tensor unit in a
symmetric monoidal closed category --- of which $\mathbf{2}$ is a
leading example.  Every
such interval $I$ gives rise to a (higher-dimensional) sesquicategory
structure on its ambient category and in some cases (such as when the
monoidal structure is cartesian) this turns out to be a 2-category
(indeed, a strict $\omega$-category).  It is our principal goal to
investigate certain properties of such an induced 2-category structure in
terms of the interval itself.  In particular, our main theorem (Theorem
\ref{theorem:rep_2cat}) gives a characterization of those strict intervals $I$ for
which the induced 2-category structure is \emph{finitely bicomplete} in
the 2-categorical sense.  A strict interval $I$ with this property is
said to be \emph{representable} and the content of Theorem
\ref{theorem:rep_2cat} is that a strict interval $I$ is representable
whenever it is a distributive lattice with top and bottom elements
which are, in a suitable sense, its generators.  

We note here that neither the closed unit
interval in the category of spaces nor the simplicial
interval in the category of simplicial sets are
examples of strict interval objects in the sense of the present
paper.  For example, although the closed unit interval can be equipped
with suitable structure maps, it fails to satisfy the defining
equations for cocategories, which are only satisfied up to homotopy.
Instead it is expected that these are examples of ``weak $\omega$-intervals''
in the sense that they are weak co-$\omega$-categories.  As such, the
present paper may be regarded as, in part, laying the groundwork for
later investigation of these intervals and the corresponding weak
higher-dimensional completeness properties of the model structures to
which they give rise.

The plan of this paper is as follows.  Section \ref{section:cocats} is
concerned with introducing the basic definitions and examples.  In
particular, we give the leading examples of strict intervals and explain the
the resulting (higher-dimensional) sesquicategory structure and when
it results in a strict 2-category.  In Section \ref{section:rep_section} we
recall the 2-categorical notion of finite bicompleteness and prove our
main results including Theorem \ref{theorem:rep_2cat}.  Lack
\cite{Lack:HA2M} has shown that every finitely bicomplete
2-category can be equipped with a model structure in which the weak
equivalences are categorical equivalences and the fibrations are
isofibrations and in Section \ref{section:isofib} we briefly explain when, given
the presence of a strict interval $I$ which is representable, this
model structure can be lifted, using a theorem due to Berger and
Moerdijk \cite{Berger:AHTO}, to the category of reduced operads.

\subsection*{Notation and conventions}

Throughout we assume, unless otherwise stated, that the ambient
category $\EE$ is a (finitely) bicomplete symmetric monoidal closed
category (for further details regarding which we refer the reader to \cite{MacLane:CWM}).  We employ common notation $(A\tensor B)$ and
$[B,A]$ for the tensor product and internal hom of objects $A$ and $B$, respectively.
We denote the tensor unit by $U$ (instead of the more common $I$) and
the natural isomorphisms associated to the symmetric monoidal closed
structure of $\EE$ are denoted by $\lambda\colon U\tensor A\to<150>A$,
$\rho\colon A\tensor U\to<150>A$, $\alpha\colon A\tensor
(B\tensor C)\to<150>(A\tensor B)\tensor C$, and $\tau\colon A\tensor
B\to<150>B\tensor A$.  Associated to the closed structure we denote
the isomorphism $[U,A]\to<150>A$ by $\partial$ and write
$\varepsilon\colon [U,A]\tensor U\to<150>A$ for the evaluation map.  We
write iterated tensor products as associating to the left so that
$A\tensor B\tensor C$ should be read as $(A\tensor B)\tensor C$.

We will frequently deal with pushouts and, if the following is a
pushout diagram
\begin{align*}
  \bfig
  \square<350,350>[C`B`A`P;f`g`g'`f']
  \efig
\end{align*}
then, when $h\colon A\to<150>X$ and $k\colon B\to<150>X$ are maps for which $h\circ
g=k\circ f$, we denote the induced map $P\to<150>X$ by $[h,k]$.
Likewise, we employ the notation $\langle h,k\rangle$ for canonical
maps into pullbacks.  

Finally, we refer the reader to \cite{Kelly:RE2C} for further details
regarding 2-categories.

\section{Intervals}\label{section:cocats}

The definition of \emph{cocategory object} in $\EE$ is exactly dual to
that of category object in $\EE$.  In order to fix notation and provide 
motivation we will rehearse the definition in full.
For us, the principal impetus for the definition of cocategories is
that a cocategory in $\EE$ provides (more than) sufficient data to
define a reasonable notion of homotopy in $\EE$ and this induced
notion of homotopy is directly related to a 2-category structure on
$\EE$.  In thinking about
cocategory objects it is often instructive to view them as analogous
to the unit interval in the category of topological spaces.  However,
the unit interval is \emph{not} a cocategory object in the category of
topological spaces and continuous functions.

\subsection{The definition}

A \myemph{cocategory} $\mathbb{C}$ in a category $\EE$ consists
objects $C_{0}$ (\myemph{object of coobjects}),  $C_{1}$
(\myemph{object of coarrows}) and $C_{2}$ (\myemph{object of
  cocomposable coarrows}) together with arrows
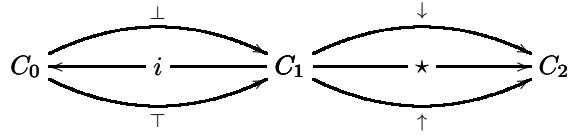
\begin{align*}
  \xy
  {\ar@/^1.25pc/^{\bottom}(0,0)*+{C_{0}};(1000,0)*+{C_{1}}};
  {\ar@/_1.25pc/_{\top}(0,0)*+{C_{0}};(1000,0)*+{C_{1}}};
  {\ar@{-}(1000,0)*+{C_{1}};(500,0)*+{i}};
  {\ar(500,0)*+{i};(0,0)*+{C_{0}}};
  {\ar@/^1.25pc/^{\downarrow}(1000,0)*+{C_{1}};(2000,0)*+{C_{2}}};
  {\ar@/_1.25pc/_{\uparrow}(1000,0)*+{C_{1}};(2000,0)*+{C_{2}}};
  {\ar@{-}(1000,0)*+{C_{1}};(1500,0)*+{\cocomp}};
  {\ar(1500,0)*+{\cocomp};(2000,0)*+{C_{2}}};
  \endxy
\end{align*}
satisfying the following list of requirements.
\begin{itemize}
\item The following square is a pushout:
  \begin{align*}
    \bfig
    \square<400,400>[C_{0}`C_{1}`C_{1}`C_{2}.;\bottom`\top`\uparrow`\downarrow]
    \efig
  \end{align*}
\item The following diagram commutes:
  \begin{align*}%\label{eq:cocat_i}
    \bfig
    \qtriangle/->`=`->/<400,350>[C_{0}`C_{1}`C_{0}.;\bottom` `i]
    \ptriangle(400,0)/<-`->`=/<400,350>[C_{1}`C_{0}`C_{0}.;\top` `]
    \efig
  \end{align*}
\item The following diagrams commute:
  \begin{align*}%\label{eq:cocat_comp}
    \bfig
    \square<400,400>[C_{0}`C_{1}`C_{1}`C_{2},;\bottom`\bottom`\cocomp`\downarrow]
    \square(1500,0)<400,400>[C_{0}`C_{1}`C_{1}`C_{2}.;\top`\top`\cocomp`\uparrow]
    \place(900,200)[\textnormal{and}]
    \efig
  \end{align*}
\item The following \myemph{co-unit laws} hold:
  \begin{align*}%\label{eq:cocat_triangle}
    \bfig
    \Atrianglepair/=`->`=`<-`->/<600,400>[C_{1}`C_{1}`C_{2}`C_{1}.;`\cocomp`
    `{[}\bottom\circ i,C_{1}{]}`{[}C_{1},\top\circ i{]}]
    \efig
  \end{align*}
\item Finally, let the object $C_{3}$ (\myemph{the object of cocomposable
    triples}) be defined as the following pushout:
  \begin{align*}
    \bfig
    \square<400,400>[C_{1}`C_{2}`C_{2}`C_{3},;\uparrow`\downarrow`\downharpoonright`\upharpoonright]
    \efig
  \end{align*}
  The \myemph{coassociative law} then states that the following diagram commutes:
  \begin{align*}%\label{eq:cocat_ass}
    \bfig
    \square<500,400>[C_{1}`C_{2}`C_{2}`C_{3}.;\cocomp`\cocomp`{[}\downharpoonright\circ\downarrow,\upharpoonright\circ\cocomp{]}`{[}\downharpoonright\circ \cocomp,\upharpoonright\circ\uparrow{]}]
    \efig
  \end{align*}
\end{itemize}
\begin{remark}
  The map $\bottom$ is the dual of the domain map, $\top$ is the dual
  of the codomain map, and $\downarrow$ and $\uparrow$
  are dual to the first and second projections, respectively.  This
  notation, and the other notation occurring in the definition, is
  justified by the interpretation
  of these arrows in the examples considered below.  We refer to $i$
  and $\cocomp$ as the \emph{coidentity} and \emph{cocomposition}
  maps, respectively.
\end{remark}
If $\mathbb{C}=(C_{0},C_{1},C_{2}) $ is a cocategory object and $A$
is any object of $\EE$, then $A\tensor\mathbb{C}=(A\tensor C_{0},A\tensor
C_{1},A\tensor C_{2})$ is also a cocategory in $\EE$.  Moreover,
if $\mathbb{C}$ is a cocategory in $\EE$ and $A$ is any object,
then one obtains a category $[\mathbb{C},A]$ in $\EE$ by taking
internal hom.

\begin{remark}
  The composites
  \begin{align*}
    [I,A]\to<350>^{t^{*}}[U,A]\to<350>^{\partial}A,
  \end{align*}
  with $t=\bottom,\top$, are denoted by $\partial_{0}$ and
  $\partial_{1}$, respectively.
\end{remark}

%%
%% COGROUPOIDS
%%
\subsection{Cocategories with additional structure}
We will be concerned with cocategories which possess
additional structure.
\begin{definition}
  \label{def:cogroupoid}
  A cocategory $\mathbb{C}$ in $\EE$ is a \myemph{cogroupoid} if
  there exists a \myemph{symmetry} or \myemph{coinverse} map $\sigma\colon C_{1}\to<150>C_{1}$ such
  that the following diagrams commute:
  \begin{align*}
    \bfig
    \Vtrianglepair/->`<-`->`->`->/<350,350>[C_{0}`C_{1}`C_{0}`C_{1};\bottom`\top`\top`\sigma`\bottom]
    \efig
  \end{align*}
  \begin{align*}
    \bfig
    \square<350,350>[C_{1}`C_{2}`C_{0}`C_{1};\cocomp`i`{[}\sigma,C_{1}{]}`\top]
    \efig
    \qquad\text{ and }\qquad
    \bfig
    \square<350,350>[C_{1}`C_{2}`C_{0}`C_{1};\cocomp`i`{[}C_{1},\sigma{]}`\bottom]
    \efig
  \end{align*}
\end{definition}
When $\mathbb{C}$ is a cogroupoid in $\EE$ and
$A$ is an object of $\EE$, $[\mathbb{C},A]$ is a groupoid in $\EE$.
\begin{definition}\label{def:interval}
  A cocategory object $\mathbb{C}$ in a category $\EE$ is
  said to be a \myemph{strict interval} if
  the object $C_{0}$ of coobjects is the tensor
  unit $U$.  When $\mathbb{C}$ is a strict interval we often write $I$
  instead of $C_{1}$ and $I_{2}$ instead of $C_{2}$.  
  When a strict interval $I$ is a cogroupoid it is said to be \myemph{invertible}.
\end{definition}
\begin{remark}
  Because we will be dealing throughout exclusively with strict
  intervals the adjective ``strict'' will henceforth be omitted.
\end{remark}
Cocategories in $\EE$ together with their obvious morphisms form a
category $\cocat(\EE)$.  There is also a category $\intervalsin(\EE)$
of strict intervals in $\EE$. 

\subsection{Examples of cocategories}\label{section:examples}

Before introducing some examples of cocategory objects it will be
useful to first record the following lemma.
\begin{lemma}\label{lemma:additive}
  Assume $\EE$ is an additive symmetric monoidal closed category.
  If we are given an object $C$ together with arrows $i\colon C\to<150>
  U$, and $\bottom,\top\colon U\to<150>C$ such that
  $i\circ\bottom=1_{U}=i\circ\top$, then there
  exists a map $\cocomp\colon  C\to<150> C_{2}$ such that the resulting structure
  is a strict interval.
  \begin{proof}
    Set $\cocomp\colon =\downarrow+\uparrow-(\uparrow\circ\bottom\circ i)$.
    The axioms for a cocategory are then immediate.
  \end{proof}
\end{lemma}
The following are examples of cocategories and intervals:
\begin{enumerate}
\item Every object $A$ of a category $\EE$ determines a cocategory
  object given by setting $A_{i}\colon =A$ for $i=0,1,2$ and defining all
  of the structure maps to be the identity $1_{A}$.  This is said to
  be the \myemph{discrete cocategory on $A$}. The discrete
  cocategory on the tensor unit $U$ is the terminal object
  in $\intervalsin(\EE)$.
\item There is an (invertible) interval in
  $\EE$ obtained by taking the object of coarrows to be $U+U$ with
  $\bottom$ and $\top$ the coproduct injections.  This is the
  initial object in $\intervalsin(\EE)$.  Indeed, a topos $\EE$ is
  Boolean if and only if its subobject classifier $\Omega$ is an
  invertible interval with $\bottom$ and
  $\top$ the usual ``truth-values''.  (Observe that if a map
  $\sigma\colon\Omega\to<150>\Omega$ in a topos satisfies
  $\sigma^{2}=1_{\Omega}$, $\sigma(\top)=\bottom$ and
  $\sigma(\bottom)=\top$, then $\sigma=\neg$.)
\item In $\cat$ the category $\mathbf{2}$ which is the free category on
  the graph consisting of two vertices and one edge between them is a
  cocategory object.  Similarly, the free groupoid $\mathbf{I}$ on
  this graph is an invertible interval in $\cat$ and in $\groupoids$
  with the following structure:
  \begin{align*}
    \bfig
    \morphism(0,0)|b|/{@{>}@/_1em/}/<400,0>[\bottom`\top;u]
    \morphism(0,0)|a|/{@{<-}@/^1em/}/<400,0>[\bottom`\top;d]
    \efig
  \end{align*}
  such that $u$ and $d$ are inverse and where
  $\bottom,\top\colon \mathbf{1}\two<150>\mathbf{I}$ are the obvious
  functors.
  $\mathbf{I}_{2}$ is then the result of gluing $\mathbf{I}$ to
  itself along the top and bottom:
  \begin{align*}
    \bfig
    \morphism(0,0)|b|/{@{>}@/_1em/}/<400,0>[\bottom`\bullet;]
    \morphism(0,0)|a|/{@{<-}@/^1em/}/<400,0>[\bottom`\bullet;]
    \morphism(400,0)|b|/{@{>}@/_1em/}/<400,0>[\bullet`\top.;]
    \morphism(400,0)|a|/{@{<-}@/^1em/}/<400,0>[\bullet`\top.;]
    \efig
  \end{align*}
  Cocomposition $\cocomp\colon \mathbf{I}\to<150>\mathbf{I}_{2}$ is the functor given by
  $\cocomp(\bottom)\colon =\bottom$ and $\cocomp(\top)\colon =\top$, and the initial and final segment
  functors are defined in the evident way.  Finally,
  $\sigma\colon \mathbf{I}\to<150>\mathbf{I}$ is defined by $\sigma(\bottom)\colon =\top$ and
  $\sigma(\top)\colon =\bottom$.  We note that these examples also
  generalize to the case of internal categories in a suitably
  complete and cocomplete category $\EE$.
\item Assume $R$ is a commutative ring (with $1$) and let $\chplus(R)$ be the category of (non-negatively graded) chain
  complexes of $R$-modules; then there exists an (invertible) interval
  $\mathbb{I}$ in $\chplus(R)$ which we now describe.  $\mathbb{I}^{0}$
  is the chain complex which consists of $R$ in degree $0$
  and is $0$ in all other degrees.  $\mathbb{I}^{1}$ is given by
  \begin{align*}
    \cdots \to<350>^{d}0\to<350>^{d}R &
    \to<350>R\oplus R\\
    x & \to/|->/<350>(x,-x),
  \end{align*}
  where $x\in R$.  $\bottom$ and $\top$ are the left and
  right inclusions, respectively.
  $i\colon \mathbb{I}^{1}\to<150>\mathbb{I}^{0}$ is given by addition in
  degree 0 and the zero map in all other degrees. Cocomposition is
  given by Lemma \ref{lemma:additive}.  The symmetry
  $\sigma\colon \mathbb{I}\to<150>\mathbb{I}$ is given by taking additive
  inverse in degree 1 and by sending $(x,y)$ to $(y,x)$ in degree 0,
  for $x,y$ in $R$.
\item Let $\rmod{R}$ denote the category of $R$-modules, for $R$ a
  commutative ring. Assume given a set $A$ together with two (not necessarily
  distinct) elements $\bottom,\top\in A$.  We obtain an interval,
  again appealing to Lemma \ref{lemma:additive}, by applying the free
  $R$-module functor $\sets\to<150>\rmod{R}$ to these data and the
  canonical map $A\to<150>1$.  
\item Following the approach from Example (5) we obtain a further
  class of examples.  Let $B$ be a bialgebra over a commutative ring $R$ with unit
  $\eta$ and counit $\epsilon$. Set $\bottom\colon =\eta$, $\top\colon =\eta$ and
  $i\colon =\epsilon$.  These data determine an interval in $\rmod{R}$ by
  Lemma \ref{lemma:additive}.
\item Assuming $\EE$ is a 2-category which is finitely cocomplete in
  the 2-categorical sense (as discussed, e.g., in Section
  \ref{section:rep_section} below) there exists for every object $A$ of $\EE$ 
  a cocategory $(A\cdot\mathbf{2})$ obtained by taking the tensor
  of $A$ with the category $\mathbf{2}$ (this fact can be
  found in its dual form in \cite{Street:FYL2}).  When $\EE$ is
  simultaneously equipped with a $\cat$-enriched symmetric monoidal
  closed structure it follows that the 2-category structure on $\EE$
  is induced, in the sense of Theorem \ref{theorem:is_2_cat}, by the
  interval $(U\cdot\mathbf{2})$.  Note that the assumption of
  $\cat$-enrichedness is necessary.
\item Consider both the cartesian and the Gray monoidal closed
  structures on the category $\textbf{2-}\cat$ of small 2-categories
  (cf.~\cite{Gordon:CT}).
  Because the tensor unit for both of these monoidal structures is the
  terminal object $\mathbf{1}$ it follows that intervals for one
  monoidal structure are the same as intervals for the other.  
  E.g., the category $\mathbf{2}$, regarded as a 2-category with no
  non-identity 2-cells, is an interval in both of these monoidal
  structures.
\end{enumerate}
As we have already noted, the topological unit interval $I=[0,1]$ in
$\topcat$ fails to satisfy the co-associativity and co-unit laws on
the nose and is therefore not an interval in the present sense.
\begin{remark}
  The question of what kinds of cocategories can exist in a topos has
  been addressed by Lumsdaine \cite{Lumsdaine:SOCC} who shows that in a
  coherent category the only cocategories are ``coequivalence
  relations''.  I.e. any such cocategory must have $\bottom$ and
  $\top$ jointly epimorphic.
\end{remark}
%%
%% HOMOTOPY
%%
\subsection{Induced sesquicategory and 2-category structures}\label{section:interval_2_cat}
The first way in which we make use of the existence of an interval
object in $\EE$ is to define homotopy.
\begin{definition}\label{def:homotopy}
  Let $I$ be an interval object in $\EE$.  A
  \myemph{homotopy (with respect to $I$)
  $\eta\colon f\twocell<150>g$ from $f$ to $g$}, for $f,g\in\EE(A,B)$, is a
  map $\eta\colon A\tensor I\to<150>B$ such that the following diagram commutes:
  \begin{align*}
    \bfig
    \square(0,0)<600,400>[A\tensor U`A\tensor
    I`A`B;A\tensor\bottom`\rho`\eta`f]
    \square(600,0)/<-`->`->`<-/<600,400>[A\tensor I`A\tensor
    U`B`A;A\tensor\top` `\rho`g]
    \efig
  \end{align*}
\end{definition}
\begin{example}\label{example:homotopy}
  Notions of homotopy corresponding to several of the intervals from
  Section \ref{section:examples} are enumerated below.
  \begin{enumerate}
  \item The discrete interval
    on $U$ generates the finest notion of homotopy (in terms of the number of homotopy
    classes of maps).  I.e., there
    exists a homotopy between $f$ and $g$ with respect to this 
    cocategory if and only if $f$ and $g$ are identical.
  \item The initial object of $\intervalsin(\EE)$ generates the
    coarsest homotopy relation: all maps are homotopic.  Indeed,
    given maps $f$ and $g$ there exists, with respect to this
    cocategory, a unique homotopy $f\twocell<150>g$.
  \item In $\cat$, homotopies $f\twocell<150>g$ with respect to
    $\mathbf{2}$ correspond to natural
    transformations $f\twocell<150>g$ and, similarly, homotopies with respect to
    $\mathbf{I}$ correspond to natural isomorphisms.
  \item In $\chplus(R)$, $\mathbb{I}$ induces the usual notion of
    chain homotopy.
  \item In the case of the interval $\mathbf{2}$ in $\textbf{2-}\cat$,
    homotopies with respect to the cartesian monoidal structure
    correspond to 2-natural transformations whereas homotopies with
    respect to the Gray monoidal structure correspond to pseudonatural
    transformations.
  \end{enumerate}
\end{example}
Recall that a \emph{sesquicategory} \cite{Street:CS} is a structure which satisfies all
of the axioms of a 2-category with the exception of the interchange
law. We will now see that, in the presence of a small amount of additional
structure, the notion of homotopy from above induces on $\EE$ the structure of a
sesquicategory.  We will also see that there are always at least two
was to endow $\EE$ with this additional structure (although there may,
as we will see, be more than these two).
\begin{definition}\label{definition:diagonal}
  Assume that $I$ is an interval in $\EE$.  We say that a map
  $\Delta\colon I\to<150>I\tensor I$ is a \myemph{diagonal (for $I$)} if the
  following conditions are satisfied:
  \begin{enumerate}
  \item $(I,\Delta,i)$ is a comonoid (cf. Appendix \ref{section:monoids_comonoids_hopf}).
  \item The diagram
    \begin{align*}
      \xy
      {\ar^{\lambda^{-1}}(0,400)*+{U};(600,400)*+{U\tensor U}};
      {\ar_{\diamond}(0,400)*+{U};(0,0)*+{I}};
      {\ar_{\Delta}(0,0)*+{I};(600,0)*+{I\tensor I}};
      {\ar^{\diamond\tensor\diamond}(600,400)*+{U\tensor
          U};(600,0)*+{I\tensor I}};
      \endxy
    \end{align*}
    commutes, for $\diamond=\bottom,\top$.
  \end{enumerate}
\end{definition}
\begin{example}\label{remark:comonoid_structure}
  Every interval $I$ in $\EE$ has an associated
  diagonal.  For the comonoid comultiplication
  $\Delta\colon I\to<150>I\tensor I$ we first form, using the
  fact that $I_{2}$ is the pushout of $\bottom$ along $\top$, the
  canonical map $[(\bottom\tensor
  I)\circ\lambda^{-1},(I\tensor\top)\circ\rho^{-1}]\colon I_{2}\to<150>I\tensor
  I$.  We then define 
  \begin{align*}
    \Delta & \colon = [(\bottom\tensor
    I)\circ\lambda^{-1},(I\tensor\top)\circ\rho^{-1}]\circ\cocomp
  \end{align*}
  With these definitions the comonoid axioms follow from the counit and coassociativity laws for
  cocategories and the second condition from Definition
  \ref{definition:diagonal} is immediate.  In the case where the
  monoidal structure on $\EE$ is cartesian, $\Delta$ is precisely
  the usual diagonal map and it is the only such map.
\end{example}
\begin{example}
  If $B$ is a bialgebra with coalgebra structure
  $(B,\Delta,\epsilon)$ as in Example (6) from Section
  \ref{section:examples}, then $\Delta$ is a diagonal.
\end{example}
\begin{example}
  If $\Delta$ is a diagonal, then so is $\tau\circ\Delta$,
  where $\tau\colon I\tensor I\to<150> I\tensor I$ is the twist map.
\end{example}
Now, assume that $I$ is an interval in $\EE$ with diagonal $\Delta$
and equip $\EE$ with the structure of a sesquicategory as follows. The
2-cells of $\EE$ are homotopies with respect to $I$.  I.e., we define
\begin{align*}
  \EE(A,B)_{1} & \colon = \EE(A\tensor I,B),
\end{align*}
which endows $\EE(A,B)$ with the structure of a category since
$[I,B]$ is an internal category in $\EE$.  Explicitly, given $\varphi$ in
$\EE(A,B)_{1}$, the domain of $\varphi$ is defined to be the arrow $\varphi\circ
(A\tensor\bottom)\circ\rho^{-1}\colon A\to<150>B$ and the codomain is $\varphi\circ
(A\tensor\top)\circ\rho^{-1}\colon A\to<150>B$.  Given arrows
$\eta\colon f\twocell<150>g$ and $\gamma\colon g\twocell<150>h$ in $\EE(A,B)$, the
vertical composite $f\twocell<150>h$ is defined as follows.  Since
$\EE$ is monoidal closed the following
square is a pushout:
\begin{align*}
  \bfig
  \square<600,400>[A\tensor U`A\tensor I`A\tensor I`A\tensor
  I_{2};A\tensor\bottom`A\tensor\top`A\tensor\uparrow`A\tensor\downarrow]
  \efig
\end{align*}
and there exists an induced map $[\eta,\gamma]\colon A\tensor
I_{2}\to<150>B$.  Recalling the third
clause from the definition of cocategory object, it is 
easily verified that $[\eta,\gamma]\circ (A\tensor \cocomp)$ is the required
vertical composite $(\gamma\cdot\eta)\colon f\twocell<150>h$.  

The horizontal composite $(\gamma*\eta)$ of a pair of 2-cells
\begin{align*}
  \bfig
  \morphism(0,0)|a|/{@{>}@/^1.5em/}/<500,0>[A`B;f]
  \morphism(0,0)|b|/{@{>}@/_1.5em/}/<500,0>[A`B;g]
  \morphism(250,75)|a|/=>/<0,-150>[`;\eta\;]
  \morphism(500,0)|a|/{@{>}@/^1.5em/}/<500,0>[B`C.;h]
  \morphism(500,0)|b|/{@{>}@/_1.5em/}/<500,0>[B`C.;k]
  \morphism(750,75)|a|/=>/<0,-150>[`;\gamma\;]
  \efig
\end{align*}
is defined to be the composite
\begin{align*}
  A\tensor I \to<250>^{A\tensor \Delta}A\tensor (I\tensor
  I)\to<150>^{\alpha}(A\tensor I)\tensor I\to<250>^{\eta\tensor
  I}B\tensor I\to<250>^{\gamma}C.
\end{align*}
The proof that the structure defined above constitutes a sesquicategory is
routine and is therefore left to the reader (the associativity and
unit laws for horizontal composition following from the coassociative
and counit laws for $\Delta$).  As such, we have the following:
\begin{proposition}\label{prop:is_sesqui_cat}
  Suppose $I$ is an interval object in $\EE$ with diagonal $\Delta$.  Then
  $\EE$ is a sesquicategory with the same objects and
  arrows, and with 2-cells the homotopies.
\end{proposition}
By virtue of Example \ref{remark:comonoid_structure}, it follows that
any interval $I$ in $\EE$ induces a sesquicategory structure on
$\EE$.  This will, however, \emph{not} in general be a 2-category
since the interchange law need not be satisfied.\footnote{The author
  is grateful to the referee for pointing out this fact and thereby
  correcting an error in the original version of this paper, and for
  suggesting the connection with part (3) of Theorem \ref{theorem:is_2_cat}.}
Proposition \ref{prop:is_sesqui_cat} has the following evident corollary:
\begin{corollary}
  \label{cor:invertible_groupoids}
  An interval $I$ in $\EE$ is invertible if and only if, for all
  objects $A$ and $B$ of $\EE$, the category $\EE(A,B)$ is a groupoid.
\end{corollary}
We will now give necessary and sufficient conditions which
characterize those cases in which the structure obtained in this way
is a genuine 2-category.
\begin{theorem}
  \label{theorem:is_2_cat}
  Assuming that $I$ is an interval in $\EE$ equipped with a diagonal
  $\Delta\colon I\to<150>I\tensor I$, then the following are equivalent:
  \begin{enumerate}
  \item $\EE$ is a 2-category, when equipped with the sesquicategory structure from
    Proposition \ref{prop:is_sesqui_cat}.
  \item The diagram
    \begin{align*}
      \begin{minipage}[h]{1.0\linewidth}
      \xy
      {\ar^{\Delta}(0,400)*+{I};(600,400)*+{I\tensor I}};
      {\ar_{\cocomp}(0,400)*+{I};(0,0)*+{I_{2}}};
      {\ar_{\Delta_{2}}(0,0)*+{I_{2}};(600,0)*+{I_{2}\tensor I_{2}}};
      {\ar^{\cocomp\tensor\cocomp}(600,400)*+{I\tensor
          I};(600,0)*+{I_{2}\tensor I_{2}}};
      \endxy  
      \end{minipage}
    \end{align*}
    commutes, where $\Delta_{2}\colon I_{2}\to<150> I_{2}\tensor I_{2}$ is the
    canonical map such that
    \begin{align*}
      \Delta_{2}\circ \downarrow & =
      (\downarrow\tensor\downarrow)\circ\Delta\\
      \Delta_{2}\circ\uparrow & = (\uparrow\tensor\uparrow)\circ\Delta.
    \end{align*}
  \item $\Delta$ is equal to the induced diagonal from Example
    \ref{remark:comonoid_structure} and $\Delta=\tau\circ\Delta$
    (i.e., $\Delta$ is cocommutative).
  \end{enumerate}
  \begin{proof}
      First, to see that (1) implies (2), assume that the interchange
      law holds to show that the diagram commutes.  Observe that we have the following diagram:
       \begin{align*}
        \xy
        {\ar@/^3em/^{\downarrow\bottom}(0,0)*+{U};(800,0)*+{I_{2}}};
        {\ar@/^3em/^{(I_{2}\tensor\downarrow\bottom)\circ\rho^{-1}}(800,0)*+{I_{2}};(1600,0)*+{I_{2}\tensor I_{2}}};
        {\ar^{\downarrow\top}(0,0)*+{U};(800,0)*+{I_{2}}};
        {\ar@/_3em/_{\uparrow\top}(0,0)*+{U};(800,0)*+{I_{2}}};
        {\ar^{(I_{2}\tensor\downarrow\top)\circ\rho^{-1}}(800,0)*+{I_{2}};(1600,0)*+{I_{2}\tensor I_{2}}};
        {\ar@/_3em/_{(I_{2}\tensor\uparrow\top)\circ\rho^{-1}}(800,0)*+{I_{2}};(1600,0)*+{I_{2}\tensor I_{2}}};
        {\ar@{=>}^{\downarrow\circ\lambda}(400,250)*+{};(400,100)*+{}};
        {\ar@{=>}^{I_{2}\tensor\downarrow}(1200,250)*+{};(1200,100)*+{}};
        {\ar@{=>}^{\uparrow\circ\lambda}(400,-75)*+{};(400,-225)*+{}};
        {\ar@{=>}^{I_{2}\tensor\uparrow}(1200,-75)*+{};(1200,-225)*+{}};
        \endxy
      \end{align*}
      and therefore the interchange law implies that 
      \begin{align*}
        \bigl((I_{2}\tensor\uparrow)*(\uparrow\circ\lambda)\bigr)\cdot\bigl((I_{2}\tensor\downarrow)*(\downarrow\circ\lambda)\bigr)
        & =
        \bigl((I_{2}\tensor\uparrow)\cdot(I_{2}\tensor\downarrow)\bigr)*\bigl((\uparrow\circ\lambda)\cdot(\downarrow\circ\lambda) \bigr)
      \end{align*}
      Since the left-hand side of this equation is equal to
      $\Delta_{2}\circ\cocomp\circ\lambda$ and, using the fact that
      $(\uparrow\cdot\downarrow)=\cocomp$, the right-hand side is
      equal to $(\cocomp\tensor\cocomp)\circ\Delta\circ\lambda$ it
      follows that the diagram in (2) commutes.

      To see that (2) implies (3), observe that the following diagram
      commutes
      \begin{align*}
        \xy
        {\ar^{\Delta}(0,400)*+{I};(600,400)*+{I\tensor I}};
        {\ar_{\cocomp}(0,400)*+{I};(0,0)*+{I_{2}}};
        {\ar_{\Delta_{2}}(0,0)*+{I_{2}};(600,0)*+{I_{2}\tensor I_{2}}};
        {\ar^{\cocomp\tensor\cocomp}(600,400)*+{I\tensor
            I};(600,0)*+{I_{2}\tensor I_{2}}};
        {\ar@{=}(600,400)*+{I\tensor I};(1200,400)*+{I\tensor I.}};
        {\ar@/_1.5em/_(.6){{[}\bottom\circ i,I{]}\tensor{[}I,\top\circ
            i{]}}(600,0)*+{I_{2}\tensor I_{2}};(1200,400)*+{I\tensor I.}};
        \endxy
      \end{align*}
      By the fact that $\Delta$ is a comonoid we therefore obtain
      \begin{align*}
        \Delta \;=\;\bigl[(\bottom\circ i\tensor
        I)\circ\Delta,(I\tensor\top\circ
        i)\circ\Delta\bigr]\circ\cocomp\;=\;\bigl[(\bottom\tensor
        I)\circ\lambda^{-1},(I\tensor\top)\circ\rho^{-1}\bigr]\circ\cocomp
      \end{align*}
      as required.  Similarly, taking $[I,\top\circ i]\tensor
      [\bottom\circ i,I]$ instead in the diagram above gives that
      $\Delta=\tau\circ\Delta$.

      To see that (3) implies (1) note that it suffices to show that,
      given a diagram
      \begin{align*}
        \xy
        {\ar@/^2em/^{f}(0,0)*+{A};(600,0)*+{B}};
        {\ar@/^2em/^{h}(600,0)*+{B};(1200,0)*+{C}};
        {\ar@/_2em/_{g}(0,0)*+{A};(600,0)*+{B}};
        {\ar@/_2em/_{k}(600,0)*+{B};(1200,0)*+{C}};
        {\ar@{=>}^{\gamma}(300,100)*+{};(300,-100)*+{}};
        {\ar@{=>}^{\delta}(900,100)*+{};(900,-100)*+{}};
        \endxy
      \end{align*}
      the two ways of composing this diagram using vertical
      composition and whiskering agree (and, in particular, agree with
      the defined horizontal composite $\delta*\gamma$).  First, one
      way of composing the diagram with whiskering and vertical
      composition is as $[\delta\circ(f\tensor
      I),k\circ\gamma]\circ(A\tensor\cocomp)$.  Using the fact that
      $f=\gamma\circ(A\tensor\bottom)\circ\rho^{-1}$ and
      $k=\delta\circ(B\tensor\top)\circ\rho^{-1}$ a straightforward
      diagram chase shows that this is equal to $\delta*\gamma$.  The
      other way of composing the diagram with whiskering and vertical
      composition gives $[h\circ\gamma,\delta\circ(g\tensor I)]\circ
      (A\tensor\cocomp)$.  Since
      $h=\delta\circ(B\tensor\bottom)\circ\rho^{-1}$ and
      $g=\gamma\circ(A\tensor\top)\circ\rho^{-1}$ a diagram chase
      gives 
      \begin{align*}
        [h\circ\gamma,\delta\circ g\tensor
        I]\circ(A\tensor\cocomp) & =\delta\circ(\gamma\tensor
        I)\circ\alpha\circ A\tensor
        [I\tensor\bottom\circ\rho^{-1},\top\tensor
        I\circ\lambda^{-1}]\circ
        A\tensor\cocomp\\
        & =\delta\circ(\gamma\tensor I)\circ\alpha\circ(A\tensor
        \tau)\circ(A\tensor\Delta).
      \end{align*}
      Therefore, since $\Delta=\tau\circ\Delta$ it follows that this
      way of composing the diagram is also equal to $\delta*\gamma$.
  \end{proof}
\end{theorem}
\begin{example}
  In the case where the monoidal structure is cartesian the equivalent
  conditions of Theorem \ref{theorem:is_2_cat} are easily seen to be
  satisfied.  I.e., in the cartesian case, $\EE$ is necessarily a 2-category.
\end{example}
\begin{example}\label{example:trivial_monoidal}
  The discrete interval (1) and the Boolean interval (2) from Section
  \ref{section:examples} both necessarily give rise to a cocommutative
  $\Delta$ and therefore also give rise to (rather trivial)
  2-categories.
\end{example}
\begin{example}
  It is easily seen that, for $\Delta$ the induced diagonal, we have
  $\tau\circ\Delta\neq\Delta$ in the category $\chplus(R)$ of chain
  complexes.  Diagrammatically, $\Delta$ and $\tau\circ\Delta$ are the maps
  \begin{align*}
    \begin{minipage}[h]{1.5cm}
      \begin{tikzpicture}
        [mynode/.style={circle,fill=black,inner sep=.5pt,minimum
          size=1mm}]
        \node[mynode] (00) at (0,0) {};
        \node[mynode] (10) at (1,0) {};
        \draw [->-] (00) to (10);
        \draw[line width=5pt, orange, opacity=.4, line
        cap=round,rounded corners=1pt] (0,0) -- (1,0);
      \end{tikzpicture}
    \end{minipage}
    \leadsto\quad
    \begin{minipage}[h]{2cm}
      \begin{tikzpicture}
        [mynode/.style={circle,fill=black,inner sep=.5pt,minimum
          size=1mm}]
        \node[mynode] (00) at (0,0) {};
        \node[mynode] (11) at (1,1) {};
        \node[mynode] (01) at (0,1) {};
        \node[mynode] (10) at (1,0) {};
        \draw [->-] (00) to (01);
        \draw [->-] (01) to (11);
        \draw [->-] (00) to (10);
        \draw [->-] (10) to (11);
        \draw[line width=5pt, orange, opacity=.4, line
        cap=round,rounded corners=1pt] (0,0) -- (0,1) -- (1,1);
      \end{tikzpicture}
    \end{minipage}
    \text{ and }\quad
    \begin{minipage}[h]{1.5cm}
      \begin{tikzpicture}
        [mynode/.style={circle,fill=black,inner sep=.5pt,minimum
          size=1mm}]
        \node[mynode] (00) at (0,0) {};
        \node[mynode] (10) at (1,0) {};
        \draw [->-] (00) to (10);
        \draw[line width=5pt, orange, opacity=.4, line
        cap=round,rounded corners=1pt] (0,0) -- (1,0);
      \end{tikzpicture}
    \end{minipage}
    \leadsto\quad
    \begin{minipage}[h]{2cm}
      \begin{tikzpicture}
        [mynode/.style={circle,fill=black,inner sep=.5pt,minimum
          size=1mm}]
        \node[mynode] (00) at (0,0) {};
        \node[mynode] (11) at (1,1) {};
        \node[mynode] (01) at (0,1) {};
        \node[mynode] (10) at (1,0) {};
        \draw [->-] (00) to (01);
        \draw [->-] (01) to (11);
        \draw [->-] (00) to (10);
        \draw [->-] (10) to (11);
        \draw[line width=5pt, orange, opacity=.4, line
        cap=round,rounded corners=1pt] (0,0) -- (1,0) -- (1,1);
      \end{tikzpicture}
    \end{minipage}
  \end{align*}
  respectively.
  As such, it follows from Theorem
  \ref{theorem:is_2_cat} that $\chplus(R)$ is
  neither a strict 2-category nor a strict $\omega$-category.
  However, there does exist an invertible 2-cell
  $(\cocomp\tensor\cocomp)\circ\Delta\iso\Delta_{2}\circ\cocomp$ which
  is given by the map
  $\varphi\colon \mathbb{I}\tensor\mathbb{I}\to<150>\mathbb{I}_{2}\tensor\mathbb{I}_{2}$
  defined as follows:
  \begin{align*}
    \varphi_{2}(a) & \colon = (0,a,0,0)\\
    \varphi_{1}(a,b,c,d) & \colon = (a+c,a,0,b,0,0,0,c,d,0,0,b+d),\text{ and}\\
    \varphi_{0}(a,b,c,d) & \colon = (a,0,b,0,c,0,0,0,d).
  \end{align*}
  Diagrammatically, this chain map is given by
  \begin{align*}
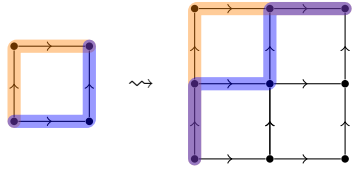
    \begin{minipage}[h]{1.5cm}
      \begin{tikzpicture}
        [mynode/.style={circle,fill=black,inner sep=.5pt,minimum
          size=1mm}]
        \node[mynode] (00) at (0,0) {};
        \node[mynode] (11) at (1,1) {};
        \node[mynode] (01) at (0,1) {};
        \node[mynode] (10) at (1,0) {};
        \draw [->-] (00) to (01);
        \draw [->-] (01) to (11);
        \draw [->-] (00) to (10);
        \draw [->-] (10) to (11);
        \draw[line width=5pt, orange, opacity=.4, line
        cap=round,rounded corners=1pt] (0,0) -- (0,1) -- (1,1);
        \draw[line width=5pt, blue, opacity=.4, line
        cap=round,rounded corners=1pt] (0,0) -- (1,0) -- (1,1);
      \end{tikzpicture}
    \end{minipage}
    \leadsto\quad
    \begin{minipage}[h]{2cm}
      \begin{tikzpicture}
        [mynode/.style={circle,fill=black,inner sep=.5pt,minimum
          size=1mm}]
        \node[mynode] (00) at (0,0) {};
        \node[mynode] (11) at (1,1) {};
        \node[mynode] (01) at (0,1) {};
        \node[mynode] (10) at (1,0) {};
        \node[mynode] (02) at (0,2) {};
        \node[mynode] (20) at (2,0) {};
        \node[mynode] (22) at (2,2) {};
        \node[mynode] (21) at (2,1) {};
        \node[mynode] (12) at (1,2) {};
        \draw [->-] (00) to (01);
        \draw [->-] (00) to (10);
        \draw [->-] (10) to (11);
        \draw [->-] (01) to (11);
        \draw [->-] (01) to (02);
        \draw [->-] (02) to (12);
        \draw [->-] (12) to (22);
        \draw [->-] (11) to (12);
        \draw [->-] (11) to (21);
        \draw [->-] (20) to (21);
        \draw [->-] (21) to (22);
        \draw [->-] (10) to (11);
        \draw [->-] (10) to (20);
        \draw[line width=5pt, orange, opacity=.4, line
        cap=round,rounded corners=1pt] (0,0) -- (0,1) -- (0,2) --
        (1,2) -- (2,2);
        \draw[line width=5pt, blue, opacity=.4, line
        cap=round,rounded corners=1pt] (0,0) -- (0,1) -- (1,1) --
        (1,2) -- (2,2);
      \end{tikzpicture}
    \end{minipage}
  \end{align*}
  where it is understood that the 2-cell of
  $\mathbb{I}\tensor\mathbb{I}$ is sent to the 2-cell in the upper
  left-hand corner of $\mathbb{I}_{2}\tensor\mathbb{I}_{2}$.
\end{example}
\begin{remark}
  As far as we know, it is an open question whether
  there exist examples, aside from the trivial ones mentioned in
  Example \ref{example:trivial_monoidal}, of intervals in the non-cartesian monoidal
  setting for which the equivalent conditions of Theorem
  \ref{theorem:is_2_cat} hold.
\end{remark}

\subsection{Semistrict higher-dimensional structure}\label{section:semistrict}

Let arrows $f,g\colon A\to<150>B$ be given and let 2-cells
$\gamma,\delta\colon f\twocell g$ also be given.  A 3-cell
$\varphi\colon \gamma\twocell \delta$ is given by an arrow
$\varphi\colon A\tensor I\tensor I\to<150> B$ which is, regarded as a 2-cell, a
homotopy $\gamma\twocell\delta$ such that 
\begin{align*}
  \xy
  {\ar^{A\tensor\bottom\tensor I}(0,400)*+{A\tensor U\tensor
      I};(1000,400)*+{A\tensor I\tensor I}};
  {\ar_{\rho\circ(\rho\tensor i)}(0,400)*+{A\tensor U\tensor
      I};(0,0)*+{A}};
  {\ar_{f}(0,0)*+{A};(1000,0)*+{B}};
  {\ar_{\varphi}(1000,400)*+{A\tensor I\tensor I};(1000,0)*+{B}};
   {\ar^{\rho\circ(\rho\tensor i)}(2000,400)*+{A\tensor U\tensor
      I};(2000,0)*+{A}};
  {\ar_{A\tensor\top\tensor I}(2000,400)*+{A\tensor U\tensor
      I};(1000,400)*+{A\tensor I\tensor I}};
  {\ar^{g}(2000,0)*+{A};(1000,0)*+{B}};
  \endxy
\end{align*}
commutes.  In general, given $n$-cells $\varphi$ and $\psi$ which are
bounded by $0$-cells $A$ and $B$, an $(n+1)$-cell $\xi\colon \varphi\twocell\psi$ is
given by an arrow $\xi\colon A \tensor I^{\tensor n}\to<150>B$ such that
\begin{align*}
  \xy
  {\ar@/^1.5em/^{A\tensor I^{\tensor k}\tensor\bottom\tensor I^{\tensor n-1-k}}(0,400)*+{A\tensor I^{\tensor
        k}\tensor U \tensor I^{\tensor n-1-k}};(1000,500)*+{A\tensor I^{\tensor n}}};
  {\ar_{\rho^{n-1-k}\circ(\rho\tensor i^{\tensor n-1-k})}(0,400)*+{A\tensor I^{\tensor
        k}\tensor U \tensor I^{\tensor n-1-k}};(300,0)*+{A\tensor I^{\tensor k}}};
  {\ar_{\partial_{0}^{n-1-k}\varphi}(300,0)*+{A\tensor I^{\tensor k}};(1000,0)*+{B}};
  {\ar_{\xi}(1000,500)*+{A\tensor I^{\tensor n}};(1000,0)*+{B}};
   {\ar^{\rho^{n-1-k}\circ(\rho\tensor i^{\tensor n-1-k})}(2000,400)*+{A\tensor I^{\tensor
        k}\tensor U\tensor I^{\tensor n-1-k}};(1700,0)*+{A\tensor I^{\tensor k}}};
  {\ar@/_1.5em/_(.4){A\tensor I^{\tensor k}\tensor\top\tensor I^{\tensor n-1-k}}(2000,400)*+{A\tensor I^{\tensor
        k}\tensor U \tensor I^{\tensor n-1-k}};(1000,500)*+{A\tensor I^{\tensor n}}};
  {\ar^{\partial_{1}^{n-1-k}\psi}(1700,0)*+{A\tensor I^{\tensor k}};(1000,0)*+{B}};
  \endxy
\end{align*}
commutes for $0\leq k\leq n-1$.  Composition of higher-dimensional cells must
be specified depending on whether the cells in question meet at a
$0$-cell or at a higher-dimensional cell. 

First, suppose given two
$(n+1)$-cells $\varphi$ and $\psi$ such that
$\partial^{n+1-k}_{1}\varphi=\partial^{n+1-k}_{0}\psi$ for $1\leq
k\leq n$.  Then the following diagram commutes
\begin{align*}
  \xy
  {\ar^(.6){A\tensor I^{\tensor k-1}\tensor\bottom\tensor I^{\tensor
        n-k}}(0,400)*+{A\tensor I^{\tensor k-1}\tensor U\tensor
      I^{\tensor n-k}};(1600,400)*+{A\tensor I^{\tensor n}}};
  {\ar_{A\tensor I^{\tensor k-1}\tensor \top\tensor I^{\tensor n-k}}(0,400)*+{A\tensor I^{\tensor k-1}\tensor U\tensor
      I^{\tensor n-k}};(0,0)*+{A\tensor I^{\tensor n}}};
  {\ar_{\varphi}(0,0)*+{A\tensor I^{\tensor n}};(1600,0)*+{B}};
  {\ar^{\psi}(1600,400)*+{A\tensor I^{\tensor n}};(1600,0)*+{B}};
  \endxy
\end{align*}
and therefore induces the map $[\varphi,\psi]\colon A\tensor I^{\tensor
  k-1}\tensor I_{2}\tensor I^{\tensor n-k}\to<150>B$.  We define the
``vertical'' composite of $\varphi$ and $\psi$ by
\begin{align*}
  \psi*_{k}\varphi & \colon = [\varphi,\psi]\circ A\tensor I^{\tensor
    k-1}\tensor\cocomp\tensor I^{\tensor n-k}.
\end{align*}
That $(\psi*_{k}\varphi)\circ ( A\tensor I^{\tensor m}\tensor \diamond
  \tensor I^{n-1-m}) = \partial_{\diamond}^{n-m}(\psi*_{k}\varphi)$,
for $\diamond=\bottom,\top$, is straightforward in the cases where
$m+1\geq k$ and is by the counit law when $m+1<k$.

Next, suppose given two $(n+1)$-cells $\varphi$ and $\psi$ such that
$\partial^{n+1}_{0}\varphi=A$,
$\partial^{n+1}_{1}\varphi=B=\partial^{n+1}_{0}\psi$, and
$\partial^{n+1}_{1}\psi=C$.  The ``horizontal'' composite of $\varphi$
and $\psi$ is given by the composite
\begin{align*}
  \xy
  {\ar^{A\tensor\Delta^{\tensor n}}(0,0)*+{A\tensor I^{\tensor
        n}};(800,0)*+{A\tensor I^{\tensor 2n}}};
  {\ar(800,0)*+{A\tensor I^{\tensor 2n}};(1600,0)*+{A\tensor
      I^{\tensor n}\tensor I^{\tensor n}}};
  {\ar^(.6){\varphi\tensor I^{\tensor n}}(1600,0)*+{A\tensor I^{\tensor
        n}\tensor I^{\tensor n}};(2500,0)*+{B\tensor I^{\tensor n}}};
  {\ar^(.6){\psi}(2500,0)*+{B\tensor I^{\tensor n}};(3000,0)*+{C}};
  \endxy
\end{align*}
where the second arrow is obtained by rearranging factors (using $\alpha$ and
$\tau$) in the obvious way so that the two $I^{\tensor n}$ in the
codomain correspond to the original $I^{\tensor n}$ in the domain of
$A\tensor \Delta^{\tensor n}$.  That $(\psi*_{0}\varphi)\circ ( A\tensor I^{\tensor m}\tensor \diamond
\tensor I^{n-1-m}) = \partial_{\diamond}^{n-m}(\psi*_{0}\varphi)$,
for $\diamond=\bottom,\top$, is straightforward.  That the associative
and unit laws are satisfied for both the vertical and horizontal
compositions just defined is essentially the same as the verification
of these laws for the 2-dimensional sesquicategory structure.
Furthermore, it is easily seen that the interchange law
\begin{align}\label{eq:interchange}
  (\psi'*_{p}\varphi')*_{q}(\psi*_{p}\varphi) & = ( \psi'*_{q}\psi)*_{p}(\varphi'*_{q}\varphi),
\end{align}
for $0<q<p\leq n$ and $(n+1)$-cells $\varphi,\varphi',\psi$ and $\psi'$,
holds.  When $q=0$ however (\ref{eq:interchange}) need not hold.
Thus, it is only with respect to this case of the interchange law that
$\EE$ fails to be a strict $\omega$-category. Thus, we have:
\begin{proposition}
  \label{prop:hd_sesquicat}
  Suppose $I$ is an interval object in $\EE$ with diagonal $\Delta$.
  Then $\EE$ is a weak higher-dimensional category in the sense that
  it satisfies all of the laws of a strict $\omega$-category except
  for (\ref{eq:interchange}) in the case $q=0$ which holds up to the
  existence of a higher-dimensional cell.
\end{proposition}
In some cases Proposition \ref{prop:hd_sesquicat} can be
strengthened.  First, we have the following two corollaries, the
proofs of which are straightforward:
\begin{corollary}\label{cor:strict_omega}
  The equivalent conditions of Theorem \ref{theorem:is_2_cat} are
  satisfied if and only if $\EE$ is a strict
  $\omega$-category.
\end{corollary}
\begin{corollary}\label{cor:iso_interchange}
  There exists an invertible 2-cell
  $(\cocomp\tensor\cocomp)\circ\Delta\iso\Delta_{2}\circ\cocomp$ if
  and only if the interchange laws (\ref{eq:interchange}), for $q=0$
  and $(n+1)$-cells $\varphi,\varphi',\psi,\psi'$, hold up to the
  existence of an invertible $(n+2)$-cell.
\end{corollary}

\section{Representability}\label{section:rep_section}

Henceforth we assume given an interval $I$ which,
together with its induced diagonal, gives rise to a 2-category (i.e.,
we assume that the equivalent conditions of Theorem
\ref{theorem:is_2_cat} are satisfied).  We now turn to the proof of our main Theorem \ref{theorem:rep_2cat}
which gives necessary and sufficient conditions under which
the 2-category structure on $\EE$ is
\emph{finitely bicomplete} in the 2-categorical sense
\cite{Street:FYL2,Gray:MMCSZ}.  We will also see that, when $\EE$ is
finitely bicomplete, $I$ can be shown to possess
additional useful structure.  For example, we will see that such an
interval is necessarily both a lattice and a \emph{Hopf object} in the
sense of Berger and Moerdijk \cite{Berger:AHTO}.  

First we recall the
2-categorical notion of finite (co)completeness due to Gray
\cite{Gray:MMCSZ} and Street \cite{Street:FYL2}.  Namely, a 2-category
$\KK$ is \myemph{finitely complete} whenever it has all finite
conical limits in the 2-categorical sense and, for each object $A$, the cotensor
$(\mathbf{2}\cotensor A)$ with the category $\mathbf{2}$ exists.
Similarly, $\KK$ is \myemph{finitely cocomplete} if and only if it
possesses all finite conical colimits and tensors $(A\cdot\mathbf{2})$
with $\mathbf{2}$ exist.  It is straightforward to verify that, when
$\EE$ possesses an interval $I$, the resulting 2-category possesses
whatever conical limits and colimits $\EE$ has in the ordinary
1-dimensional sense:
\begin{lemma}\label{lemma:2limits}
  Assume that $\EE$ possesses an interval $I$ and regard $\EE$ as a
  2-category with respect to the 2-category structure induced by $I$.
  Then the conical (co)limit of a functor $F\colon \CC\to<150>\EE$ from a (small) category $\CC$
  exists if and only if the ordinary 1-dimensional (co)limit of $F$ exists.
\end{lemma}
In order to show that the 2-category structure on $\EE$
induced by an interval $I$ is finitely bicomplete it suffices, by Lemma
\ref{lemma:2limits}, to prove that tensor and cotensor products with the category
$\mathbf{2}$ exist.  Indeed, if $(\mathbf{2}\cotensor A)$ exists, it
is necessarily isomorphic to the internal hom $[I,A]$ since the 
2-natural isomorphism
\begin{align}\label{eq:cotensor_iso}
  \EE(B,\mathbf{2}\cotensor A) & \iso \EE(B,A)^{\mathbf{2}}
\end{align}
of categories restricts to a natural isomorphism of their
respective collections of objects:
\begin{align*}
  \EE(B,\mathbf{2}\cotensor A) & \iso \EE(B\tensor I,A).
\end{align*}
Similar reasoning implies that when the tensor product
$(A\cdot\mathbf{2})$ exists it is necessarily $(A\tensor I)$.
Note though that it does not \emph{a priori} follow that $[I,A]$ is
$(\mathbf{2}\cotensor A)$ in the sense of possessing the full
2-categorical universal property of $(\mathbf{2}\cotensor A)$, and
similarly for $(A\tensor I)$ and $(A\cdot\mathbf{2})$.  This remark
should be compared with the familiar fact that a 2-category with all
1-dimensional conical limits need not possesses all 2-dimensional
conical limits (cf.~\cite{Kelly:BCECT}).

As the reader may easily verify, if $I$ is an interval in $\EE$, then
there exist isomorphisms of categories
\begin{align*}
  \EE(B\tensor I,A) & \iso \EE\bigl(B,[I,A]\bigr)
\end{align*}
natural in $A$ and $B$.  Thus, it follows that $\EE$ possesses tensors
with $\mathbf{2}$ if and only if it possesses cotensors with
$\mathbf{2}$.  
\begin{definition}\label{def:rep_interval}
  An interval $I$ in $\EE$ is \myemph{representable}
  if cotensors with $\mathbf{2}$ exist with respect to the 2-category
  structure on $\EE$ induced by $I$.
\end{definition}
Thus, an interval $I$ is representable if and only if $\EE$ is a
finitely bicomplete 2-category with respect to the induced
2-category structure of Section \ref{section:interval_2_cat}.  In
particular, when $I$ is representable the monoid structure is
$\cat$-monoidal and $I$ is necessarily obtained as in Example (7) from
Section \ref{section:examples}.

\subsection{Injective boundaries}\label{section:double}

An arrow $\varphi\colon B\tensor I\tensor I\to A$ in $\EE$ determines a square
\begin{align*}
  \bfig
  \square<400,400>[\varphi_{00}`\varphi_{10}`\varphi_{01}`\varphi_{11};\varphi_{\bullet
    0}`\varphi_{0\bullet}`\varphi_{1\bullet}`\varphi_{\bullet 1}]
  \efig
\end{align*}
in $\EE(B,A)$, where our notation should be clear (e.g., $\varphi_{0\bullet}$ is
the result of precomposing $\varphi$ with $(B\tensor\bottom\tensor
I)\circ(\rho^{-1}\tensor I)$).  Because $\EE$ is assumed to satisfy the equivalent
conditions of Theorem \ref{theorem:is_2_cat} it follows that this
diagram commutes.  Consequently, there exists an induced homomorphism 
\begin{align*}
  \EE(B,A)^{\flat}\to^{\partial}\EE(B,A)^{\natural}
\end{align*}
of double categories where $\EE(B,A)^{\natural}$ is the double category
with objects the objects of $\EE(B,A)$, vertical and horizontal
arrows both the arrows in $\EE(B,A)$, and 2-cells commutative diagrams
in $\EE(B,A)$; and where $\EE(B,A)^{\flat}$ has the same objects,
horizontal and vertical arrows as $\EE(B,A)^{\natural}$, but with
2-cells arrows $\varphi\colon B\tensor I\tensor I\to A$.  Note that horizontal
composition of composable 2-cells $\varphi,\psi$ in
$\EE(B,A)^{\flat}$ is given by
$\psi\circ_{h}\varphi\colon =[\varphi,\psi]\circ(B\tensor\cocomp\tensor I)$
and vertical composition is given by
$\psi\circ_{v}\varphi\colon =[\varphi,\psi]\circ(B\tensor
I\tensor\cocomp)$.  Let $P$ be the pushout of $\cocomp\colon
I\to<150>I_{2}$ along itself.  Then maps $B\tensor P\to<150>A$ are in
bijective correspondence with the 2-cells of $\EE(B,A)^{\natural}$.
The maps $\Delta,\tau\circ\Delta\colon I\to<150>I\tensor I$ induce, by
their definitions, a canonical map $P\to<150>I\tensor I$ and the
action of $\partial$ on 2-cells is induced by precomposition with this
map.

Moving from double categories to categories, $\partial$ restricts to a functor 
\begin{align*}
  \EE(B,[I,A]) \to^{\Phi} \EE(B,A)^{\mathbf{2}}
\end{align*}
which acts by transpose under the tensor-hom adjunction.  I.e., given an
object $f$ of $\EE(B,[I,A])$, the arrow $\Phi(f)$ in $\EE(B,A)$ is
defined to be the transpose $\tilde{f}\colon B\tensor I\to<150>A$ of $f$.
Similarly, for an arrow $\varphi\colon f\twocell<150>g$ in $\EE(B,[I,A])$, $\Phi(f)$ is
obtained by projecting the transpose $\tilde{\varphi}$ to the
commutative square $\partial(\tilde{\varphi})$:
\begin{align}\label{eq:bdry_diagram}
  \bfig
  \square<400,400>[\partial_{0}f`\partial_{0}g`\partial_{1}f`\partial_{1}g.;\partial_{0}\circ\varphi`\tilde{f}`\tilde{g}`\partial_{1}\circ\varphi]
  \efig
\end{align}
The following lemma implies that if $I$ is representable, then $\Phi$ is necessarily the
natural isomorphism witnessing this fact.
\begin{lemma}
  \label{lemma:domain_codomain}
  If $I$ is representable, then, for all objects $A$ and $B$ of $\EE$,
  the functors $\Phi$ give isomorphisms of categories which are
  natural in $A$ and $B$.  Furthermore, the following diagram in $\cat$ commutes:
  \begin{align*}
    \bfig
    \Vtriangle<400,350>[\EE\bigl(B,{[}I,A{]}\bigr)`\EE(B,A)^{\mathbf{2}}`\EE(B,A);\Phi`\EE(B,\partial_{i})`\partial_{i}]
    \efig
  \end{align*}
  when $i=0,1$.
\end{lemma}
In particular, representability of $I$ is equivalent to $\Phi$ being an isomorphism
of categories natural in $A$ and $B$, which is equivalent to
$\partial$ being an isomorphism of double categories which is
similarly natural in $A$ and $B$.  Because naturality is immediate by
definition all of these are equivalent to the canonical map
$P\to<150>I\tensor I$ being an isomorphism.
\begin{definition}
  An interval $I$ has \myemph{injective boundaries} if
  2-cells $\varphi$ in the double categories $\EE(B,A)^{\flat}$ are completely determined by their
  boundaries.  I.e., for any objects $A$ and $B$ of $\EE$ and any
  2-cells $\varphi$ and $\psi$ in $\EE(B,A)^{\flat}$,
  $\partial(\varphi)=\partial(\psi)$ implies that $\varphi=\psi$.
\end{definition}
It is worth remarking that it suffices to test maps $I\tensor I\to A$
in order to determine whether or not $I$ has injective boundaries.  In
more homotopic terms, $I$ has injective boundaries provided that for
all paths $f$ and $g$ from $a$ to $b$ in $A$ there is
at most one homotopy rel endpoints $f\simeq g$.  The following
observation is a trivial consequence of the discussion above:
\begin{lemma}
  \label{lemma:boundaries}
  All representable intervals $I$ have injective boundaries.
\end{lemma}

\subsection{Lattice structure of representable intervals}\label{section:lattice_and_hopf}

We will now prove that if $I$ is representable, then it is necessarily
a unital distributive lattice in the sense of Appendix
\ref{section:monoids_comonoids_hopf}.
\begin{proposition}\label{prop:rep_meet_join}
  If $I$ is representable, then it possesses the
  structure of a unital distributive lattice
  such that $\bottom$ is the unit for join $\vee\colon I\tensor I\to<150>I$
  and $\top$ is the unit for meet $\wedge\colon I\tensor I\to<150>I$.
  Moreover, this structure is unique in the strong sense that
  meet and join are the canonical maps $I\tensor I\to<150>I$ such
  that both $\wedge$ and $\wedge\circ\tau$
  are 2-cells $\bottom\circ i\twocell<150>1_{I}$, and both $\vee$ and
  $\vee\circ\tau$ are 2-cells $1_{I}\twocell<150>\top\circ i$.
  \begin{proof}
    Because $I$ is representable it follows that there exists a 2-natural isomorphism
    \begin{align}\label{eq:special_case_cotensor_iso}
      \EE(U,I)^{\mathbf{2}}\to<350>^{\iso}\EE\bigl(U,[I,I]\bigr)
    \end{align}
    of categories which is given at the level of objects by
    exponential transpose.  In $\EE(U,I)$ the following diagram commutes
    \begin{align*}
      \bfig
      \square<400,400>[\bottom`\bottom`\bottom`\top;\bottom\circ
      i\circ\lambda`\bottom\circ i\circ\lambda`\lambda`\lambda]
      \efig
    \end{align*}
    Thus, applying (\ref{eq:special_case_cotensor_iso}) to this arrow of
    $\EE\bigl(U,I)^{\mathbf{2}}$ yields a map $\boxdot\colon U\tensor
    I\to<150>[I,I]$.  Denote by $\wedge\colon I\tensor I\to<150>I$ the
    transpose of the composite $\boxdot\circ\lambda^{-1}$ and observe, by
    definition and Lemma \ref{lemma:domain_codomain}, that $\top$ is a
    unit for this operation and that the diagram
    \begin{align}\label{eq:meet_condition_2}
      \bfig
      \square(0,0)<600,400>[U\tensor I`I\tensor I`I`I;\bottom\tensor
      I`\lambda`\wedge`\bottom\circ i]
      \square(600,0)/<-`->`->`<-/<600,400>[I\tensor I`I\tensor
      U`I`I;I\tensor \bottom``\lambda`\bottom\circ i]
      \efig
    \end{align}
    also commutes.  In the same way, applying the isomorphism
    (\ref{eq:special_case_cotensor_iso}) to the arrow
    \begin{align*}
      \bfig
      \square<400,400>[\bottom`\top`\top`\top;\lambda`\lambda`\top\circ
      i\circ\lambda`\top\circ i\circ\lambda]
      \efig
    \end{align*}
    of $\EE(U,I)^{\mathbf{2}}$ yields a map
    $\boxplus\colon U\tensor I\to<150>[I,I]$ for which the transpose
    $\vee\colon I\tensor I\to<150>I$ of $\boxplus\circ\lambda^{-1}$
    is an operation which has as a unit $\bottom$ and satisfies the dual
    of (\ref{eq:meet_condition_2}).  Moreover, by Lemma
    \ref{lemma:boundaries}, it follows that $\vee$ and $\wedge$ are
    the canonical maps $I\tensor I\to<150>I$ with these properties.
    For example, the idempotent law which
    states that $\vee\circ\Delta=1_{I}$ holds since
    \begin{align*}
      \vee\circ\Delta & = [\vee\circ(\bottom\tensor
      I)\circ\lambda^{-1},\vee\circ(I\tensor\top)\circ\rho^{-1}]\circ\cocomp\\
      & = [1_{I},\top\circ i]\circ\cocomp\\
      & = 1_{I}
    \end{align*}
    where the final equation is by the cocategory counit law.  The
    other idempotent law is similar.  Commutativity of the additional
    diagrams for distributive lattices also follow from Lemma
    \ref{lemma:boundaries} by a routine (but lengthy) series of
    diagram chases.
  \end{proof}
\end{proposition}
Using join $\vee\colon I\tensor I\to<150>I$ we see that $I$ is a commutative
Hopf object in the sense of \cite{Berger:AHTO} (see Appendix
\ref{section:monoids_comonoids_hopf} for the definition).
\begin{corollary}\label{cor:hopf}
  If $I$ is representable, then it is a commutative Hopf object.
  \begin{proof}
    As we have already seen $(I,\Delta,i)$ is a comonoid and both
    $(I,\vee,\bottom)$ and $(I,\wedge,\top)$ are commutative monoids.  In
    fact, $I$ can be made into a commutative Hopf object using either
    of these monoid structures.  To see this it remains to verify that
    $\vee$ and $\bottom$, as well as $\wedge$ and $\top$, are comonoid
    homomorphisms.  Since
    $\cocomp\circ\bottom=\downarrow\circ\bottom$ it follows that
    $\bottom$ is a homomorphism.  $\vee$ is seen to be a homomorphism
    by testing on boundaries.  A dual proof shows that $\wedge$ and
    $\top$ are also comonoid homomorphisms.
  \end{proof}
\end{corollary}
\begin{remark}
  We note that if $\varphi\colon I\to<150>H$ is an arrow in
  $\intervalsin(\EE)$ between representable intervals,
  then it is necessarily also a morphism of Hopf objects provided that
  $H$ and $I$ are both equipped with ``meet'' (respectively, ``join'') Hopf object structures
  from the proof of Corollary \ref{cor:hopf}.
\end{remark}

\subsection{The characterization of representable intervals}

We would now like to investigate the extent to which Proposition
\ref{prop:rep_meet_join} characterizes representable intervals.  For
the remainder of this section, unless otherwise stated we do
\emph{not} assume that $I$ is representable.  We do however assume
that there exist meet $\wedge\colon I\tensor I\to<150>I$ and join $\vee\colon I\tensor I\to<150>I$
operations which have $\bottom$ and $\top$ as respective units and
satisfy condition (\ref{eq:unital_lattice}) from Appendix
\ref{section:monoids_comonoids_hopf} (equivalently, both $\wedge$ and
$\wedge\circ\tau$ are 2-cells $\bottom\circ i\twocell<150>1_{I}$, and
both $\vee$ and $\vee\circ\tau$ are 2-cells $1_{I}\twocell<150>\top\circ i$).  

Let us recall some double category machinery from
\cite{Brown:DC2CTSC}.  An double category $\DD$ is
\myemph{edge symmetric} when it has the same horizontal and vertical
edges.  
\begin{definition}\label{def:connection}
  A \myemph{connection} on an edge symmetric double category $\DD$
  consists of maps
  \begin{align*}
    \Gamma,\Gamma'\colon \DD_{1}\to<150>\DD_{2}
  \end{align*}
  such that $\Gamma(f)$ and $\Gamma'(f)$ have boundaries
  \begin{align*}
    \xy
    {\ar^{f}(0,400)*+{\partial_{0}(f)};(500,400)*+{\partial_{1}(f)}};
    {\ar_{f}(0,400)*+{\partial_{0}(f)};(0,0)*+{\partial_{1}(f)}};
    {\ar^{1_{\partial_{1}(f)}}(500,400)*+{\partial_{1}(f)};(500,0)*+{\partial_{1}(f)}};
    {\ar_{1_{\partial_{1}(f)}}(0,0)*+{\partial_{1}(f)};(500,0)*+{\partial_{1}(f)}};
    {(900,250)*+{\text{and}}};
    {\ar^{1_{\partial_{0}(f)}}(1300,400)*+{\partial_{0}(f)};(1800,400)*+{\partial_{0}(f)}};
    {\ar_{1_{\partial_{0}(f)}}(1300,400)*+{\partial_{0}(f)};(1300,0)*+{\partial_{0}(f)}};
    {\ar^{f}(1800,400)*+{\partial_{0}(f)};(1800,0)*+{\partial_{1}(f)}};
    {\ar_{f}(1300,0)*+{\partial_{0}(f)};(1800,0)*+{\partial_{1}(f)}};
    \endxy
  \end{align*}
  respectively, for $a$ in $\DD_{1}$, and such that $\Gamma$ and
  $\Gamma'$ satisfy several further conditions which we now describe.
  First, they are required to preserve identities in the sense that
  $\Gamma(1_{a})=1_{1_{a}}=\Gamma'(1_{a})$ for $a$ an object of
  $\DD$.  Next, it is required that, for arrows $f\colon a\to<150>b$ and
  $g\colon b\to<150>c$, $\Gamma(g\circ f)$ and $\Gamma'(g\circ f)$ are equal to the
  composites
   \begin{align*}
     \begin{minipage}{1.0\linewidth}
      \xy
      {\ar^{f}(0,800)*+{a};(500,800)*+{b}};
      {\ar^{g}(500,800)*+{b};(1000,800)*+{c}};
      {\ar_{f}(0,800)*+{a};(0,400)*+{b}};
      {\ar_{1_{b}}(0,400)*+{b};(500,400)*+{b}};
      {\ar^{1_{b}}(500,800)*+{b};(500,400)*+{b}};
      {\ar_{g}(500,400)*+{b};(1000,400)*+{c}};
      {\ar^{1_{c}}(1000,800)*+{c};(1000,400)*+{c}};
      {\ar_{g}(0,400)*+{b};(0,0)*+{c}};
      {\ar_{1_{c}}(0,0)*+{c};(500,0)*+{c}};
      {\ar_{1_{c}}(500,0)*+{c};(1000,0)*+{c}};
      {\ar^{1_{c}}(1000,400)*+{c};(1000,0)*+{c}};
      {\ar^{g}(500,400)*+{b};(500,0)*+{c}};
      {\ar_{g}(500,400)*+{b};(1000,400)*+{c}};
      {(250,600)*+{\Gamma(f)}};
      {(750,600)*+{1^{v}_{g}}};
      {(250,200)*+{1^{h}_{g}}};
      {(750,200)*+{\Gamma(g)}};
      \endxy
    \end{minipage}
    \quad\text{and}\quad
   \begin{minipage}{1.0\linewidth}
      \xy
      {\ar^{1_{a}}(0,800)*+{a};(500,800)*+{a}};
      {\ar^{1_{a}}(500,800)*+{a};(1000,800)*+{a}};
      {\ar_{1_{a}}(0,800)*+{a};(0,400)*+{a}};
      {\ar_{f}(0,400)*+{a};(500,400)*+{b}};
      {\ar^{f}(500,800)*+{a};(500,400)*+{b}};
      {\ar_{1_{b}}(500,400)*+{b};(1000,400)*+{b}};
      {\ar^{f}(1000,800)*+{a};(1000,400)*+{b}};
      {\ar_{1_{a}}(0,400)*+{a};(0,0)*+{a}};
      {\ar_{f}(0,0)*+{a};(500,0)*+{b}};
      {\ar_{g}(500,0)*+{b};(1000,0)*+{c}};
      {\ar^{g}(1000,400)*+{b};(1000,0)*+{c}};
      {\ar^{1_{b}}(500,400)*+{b};(500,0)*+{b}};
      {\ar_{1_{b}}(500,400)*+{b};(1000,400)*+{b}};
      {(250,600)*+{\Gamma'(f)}};
      {(750,600)*+{1^{h}_{g}}};
      {(250,200)*+{1^{v}_{g}}};
      {(750,200)*+{\Gamma'(g)}};
      \endxy
    \end{minipage}
  \end{align*}  
  respectively.  Finally, we require that $\Gamma$ and $\Gamma'$ are
  inverse to one another in the sense that
  $\Gamma(f)\circ_{h}\Gamma'(f)=1_{f}^{v}$ and $\Gamma'(f)\circ_{v}\Gamma(f)=1^{h}_{f}$.
\end{definition}
We will make use of the following result in the proof of our main theorem.
\begin{theorem}[Brown and Mosa (Corollary 4.4 in \cite{Brown:DC2CTSC})]\label{theorem:BandM}
  On an edge symmetric double category $\DD$, connections correspond
  to morphisms $\Theta\colon\Box\DD\to<150>\DD$ which are
  identity on objects and arrows, where $\Box\DD$ is the
  double category which has the same objects and 1-cells as $\DD$, but
  with 2-cells given by commutative squares.
\end{theorem}
We now turn to our main theorem.
\begin{theorem}\label{theorem:rep_2cat}
  An interval $I$ in $\EE$ is representable if and only if it has
  injective boundaries and possesses binary meet and join operations
  such that both $\wedge$ and $\wedge\circ\tau$ are 2-cells $\bottom\circ i\twocell<150>1_{I}$,
  and both $\vee$ and $\vee\circ\tau$ are 2-cells
  $1_{I}\twocell<150>\top\circ i$.
  \begin{proof}
    It follows from Proposition
    \ref{prop:rep_meet_join} and Lemma \ref{lemma:boundaries} that a
    representable interval possesses the required properties.  

    For the other direction of the equivalence it suffices, by the discussion in
    Section \ref{section:double}, to prove that
    $\partial\colon\EE(B,A)^{\flat}\to<150>\EE(B,A)^{\natural}$, which is immediately
    seen to be natural in $A$ and $B$, is an isomorphism of
    double categories.  For this, we first observe that using $\vee$ and
    $\wedge$ we obtain a connection on the double category $\EE(B,A)^{\flat}$ by
    letting $\Gamma$ and $\Gamma'$ send $f\colon B\tensor I\to<150>A$ to
    the composites $f\circ(B\tensor\vee)\circ\alpha^{-1}$ and
    $f\circ(B\tensor\wedge)\circ\alpha^{-1}$, respectively.  By
    definition, $\Gamma(f)$ and $\Gamma'(f)$ have the correct
    boundaries and, because $I$ has injective boundaries, the
    remaining conditions on a connection are also satisfied.  By
    Theorem \ref{theorem:BandM}, there exists a map $\Theta$ of
    double categories
    $\Box\EE(B,A)^{\flat}=\EE(B,A)^{\natural}\to\EE(B,A)^{\flat}$
    which is identity on objects and arrows.  In particular, $\Theta$
    is a section of $\partial$ and therefore, by injective boundaries,
    $\partial$ is an isomorphism of double categories.
  \end{proof}
\end{theorem}
Although most of the examples of intervals studied earlier are already
known to give rise to finitely bicomplete 2-category structures, it is
nonetheless instructive to consider these cases in light of the theorem.
\begin{example}
  Consider the following intervals:
  \begin{enumerate}
  \item The interval $I$ obtained by taking the discrete cocategory
    on the tensor unit $U$ is representable, with meet and join both
    the structure map $\lambda=\rho\colon U\tensor U\to<150>U$.
  \item Using the isomorphism $(U+U)\tensor(U+U)\iso(U+U)+(U+U)$ it is
    easily seen that the interval $(U+U)$ satisfies the necessary and
    sufficient conditions from Theorem \ref{theorem:rep_2cat} for
    being representable and therefore gives rise to a finitely
    bicomplete 2-category.
  \item In $\cat$ both $\mathbf{2}$ and $\mathbf{I}$ are
    representable.  Of course, this can be easily verified directly,
    but one can also check that the hypotheses of the theorem are
    satisfied.  For instance, in both cases the meet map $\wedge$ is
    the functor which sends an object $(s,t)$ of $\mathbf{2}\times\mathbf{2}$ to
    $\top$ if $s=t=\top$ and to $\bottom$ otherwise.
  \item  We will now give an example of an interval giving rise to a
    2-category structure, but which is \emph{not} representable.  Let us work in $\cat$ with the cartesian monoidal structure.  Let $L$
    be the free category on the graph with one vertex $\mu$ and one edge
    $\omega\colon \mu\to\mu$.  (I.e., it is the free monoid on a single
    generator.) Then $L_{2}$ is the free category with one vertex $\mu$
    and two $l,r\colon \mu\to\mu$ and $\downarrow(\omega)=l$,
    $\uparrow(\omega)=r$, $\cocomp(\omega)=r\circ l$.  This is an interval
    in $\cat$ and it induces a trivial notion of homotopy.  Namely, for
    functors $F,G\colon A\to B$ if $F\simeq G$, then $F=G$ and there exists
    for each $a$ in $A$ a loop $\varphi_{a}\colon Fa\to Fa$ such that
    \begin{align*}
      \varphi_{b}\circ Ff & = Ff\circ\varphi_{a}
    \end{align*}
    for $f\colon a\to b$ in $A$.  Roughly, this interval generates the same
    notion of homotopy as the discrete interval, \emph{but} the data of a
    homotopy for $L$ is not the same as the data of a homotopy for the
    discrete interval.  As such, the resulting 2-category structures are
    not (\emph{a priori}) the same.
    
    The interval $L$ can possess neither meet nor join operations since
    in this case $\bottom=\top$ and so we would have
    \begin{align*}
      1_{\mu}=\bottom=\omega\wedge\bottom=\omega\wedge\top=\omega
    \end{align*}
    which is false.  Thus, by the characterization theorem it follows that
    $L$ is not representable.
  \end{enumerate}
\end{example}

\section{Homotopy theoretic consequences}\label{section:isofib}

The purpose of this section is to relate the considerations on
intervals from the foregoing sections to several known results
from homotopy theory.  In particular, we show that, under suitable
hypotheses on $\EE$, if $I$ is a representable interval in $\EE$, then
the  ``isofibration'' model structure on $\EE$ due to Lack
\cite{Lack:HA2M} can be lifted to the category of (reduced) operads
using a theorem of Berger and Moerdijk \cite{Berger:AHTO}.  In
order to apply the machinery of \emph{ibid} it is first necessary to
construct a \emph{Hopf interval}, which is essentially a cylinder
object equipped with the structure of a Hopf object.  As such, the
principal observation in this section is that, when $\EE$ is
cocomplete in the 1-dimensional sense, it is possible to construct the
free Hopf object generated by the interval $I$.  We refer the reader
to \cite{Hovey:MC} for more information regarding model categories.

Although we will not consider those intervals $I$ which fail to be
representable (or to give rise to 2-categories) in our discussion of
homotopy theory below, we would
like to mention that some effort has been made to investigate the
homotopy theory of intervals arising in the setting of such categories
as the category of chain complexes.  In particular, Stanculescu
\cite{Stanculescu:HTESM} has employed intervals in his work on the
homotopy theory of categories enriched in simplicial modules.

\subsection{The isofibration model structure}

Now, assuming (as we will throughout the remainder of this
section) that $\EE$ is a finitely bicomplete symmetric monoidal closed
category with a representable interval $I$, it follows from a theorem
due to Lack \cite{Lack:HA2M} that $\EE$ can be equipped with a model
structure in which the weak equivalences are the categorical
equivalences and the fibrations are isofibrations.
Recall that an arrow $f\colon A\to<150>B$ in a 2-category is said to be a
\myemph{categorical equivalence} if there exists a map $f'\colon B\to<150>A$
together with isomorphisms $f\circ f'\iso 1_{B}$ and $f'\circ f\iso
1_{A}$.  A functor $F\colon \CC\to<150>\DD$ in $\cat$ is said to be an isofibration when
isomorphisms in $\DD$ whose codomains lie in the image of $F$ can be
lifted to isomorphisms in $\CC$.  This notion also makes sense in
arbitrary 2-categories $\EE$.  We define a map
$f\colon A\to<150>B$ in $\EE$ to be an
\myemph{isofibration} if, for any object $E$ of $\EE$, the induced map
\begin{align*}
  \EE(E,A)\to<350>^{f_{*}}\EE(E,B)
\end{align*}
is an isofibration in $\cat$.  
\begin{definition}\label{def:model_cat_cat}
  Assume $\EE$ is a finitely bicomplete 2-category with a model
  structure.  Then $\EE$ is a \myemph{model $\cat$-category} provided
  that if $p\colon E\to<150>B$ is a fibration and $i\colon X\to<150>Y$ is a
  cofibration, then the induced functor 
  \begin{align*}
    \EE(Y,E)\to<350>^{\langle p_{*},i^{*}\rangle}\EE(Y,B)\times_{\EE(X,B)}\EE(X,E)
  \end{align*}
  is an isofibration which is simultaneously an equivalence if either
  $p$ or $i$ is a weak equivalence.
\end{definition}
With these definitions, Lack
\cite{Lack:HA2M} proved the following theorem:
\begin{theorem}[Lack]
  \label{theorem:Lack}
  If $\EE$ is a finitely bicomplete 2-category, then it bears the
  structure of a model
  $\cat$-category in which the weak equivalences are the
  equivalences, the fibrations are the isofibrations and the
  cofibrations are those maps having the left-lifting property with
  respect to maps which are simultaneously fibrations and weak
  equivalences.
\end{theorem}
We will refer to such a model structure on a 2-category $\EE$ as the
\myemph{isofibration model structure} on $\EE$.  Every object is both fibrant
and cofibrant in this model structure.  It is an immediate
consequence of Theorem \ref{theorem:Lack} that
when $\EE$ is a finitely bicomplete symmetric monoidal closed category with a representable
interval $I$ it is also a model $\cat$-category with
the isofibration model structure.

\subsection{The free Hopf interval generated by $I$}

In \cite{Berger:AHTO}, a (commutative) \myemph{Hopf interval} in a symmetric
monoidal model category is defined to be a cylinder object 
\begin{align*}
  \bfig
  \Vtriangle<300,300>[U+U`H`U;`\nabla`]
  \efig
\end{align*}
on $U$ such that $H$ is a (commutative) Hopf object, and both maps $U+U\to<150>H$ and $H\to<150>U$ are
homomorphisms of Hopf objects, where to be a cylinder object means
that $U+U\to<150>H$ is a cofibration and $H\to<150>U$ is a weak
equivalence. Here $U$ has the trivial Hopf object
structure given by by structure map $\lambda\colon U\tensor U\to<150>U$ and
its inverse.  On the other hand, $(U+U)$ is given the structure of
a commutative Hopf object using the ``meet'' Hopf object structure
described in Corollary \ref{cor:hopf} (cf.~Example (2) from Section
\ref{section:examples}), which coincides with the Hopf
object structure on described in \emph{ibid}.  We emphasize
that a Hopf interval need not be an interval in the sense of
Definition \ref{def:interval}.  The following lemma shows that we
cannot in general expect $I$ itself to be a Hopf interval in the
isofibration model structure.
\begin{lemma}\label{lemma:i_equiv}
  The following are equivalent:
  \begin{enumerate}
  \item $I$ is invertible.
  \item The structure map $\lambda\colon U\tensor I\to<150>I$, regarded as a
    2-cell $\bottom\twocell<150>\top$, possesses an inverse
    $\neg\colon \top\twocell<150>\bottom$.
  \item The meet operation, regarded as a 2-cell $\bottom\circ
    i\twocell<150>1_{I}$, has an inverse. (Or, dually, the join
    operation has an inverse.)
  \item The map $i\colon I\to<150>U$ is an equivalence.
  \end{enumerate}
  \begin{proof}
    (1) and (2) are equivalent since $\neg$ is defined using the
    existence of a coinverse map $\sigma\colon I\to<150>I$ to be
    $\sigma\circ\lambda$, and, going the other way, $\sigma$ is
    defined in terms of $\neg$ as $\neg\circ\lambda^{-1}$.  To see
    that (2) implies (3), define $\wedge'\colon I\tensor I\to<150>I$ to be
    the composite
    \begin{align*}
      I\tensor I\to<350>^{I\tensor\lambda^{-1}}I\tensor(U\tensor
      I)\to<350>^{I\tensor\neg}I\tensor I\to<350>^{\wedge}I.
    \end{align*}
    It is then easily seen that $\wedge'$ is the inverse of $\wedge$
    as arrows in the category $\EE(I,I)$.  

    That (3) implies (4) follows from the fact that $\wedge$ is a
    2-cell $\bottom\circ i\twocell<150>1_{I}$ and therefore the
    existence of an inverse for this homotopy implies that $i$ is an
    equivalence.  Going the other way, to see that (4) implies (2),
    assume given $k\colon U\to<150>I$ together with isomorphisms
    $\varphi\colon k\circ i\twocell<150>1_{I}$ and
    $\psi\colon 1_{U}\twocell<150>i\circ k$.  Then we define the inverse
    $\neg\colon \top\twocell<150>\bottom$ of $\lambda$ in $\EE(U,I)$ as
    follows.  We first form the vertical composite
    $\xi\colon =(\varphi*\bottom)$ which is, by definition, a 2-cell
    $k\twocell<150>\bottom$.  Similarly, we define
    $\zeta\colon \top\twocell<150>k$ to be the vertical composite
    $(\varphi^{-1}*\top)$.  We then set $\neg\colon =(\xi\cdot\zeta)$.
    That $\neg$ is the inverse of $\lambda$ is seen to hold by
    straightforward calculations.  For example, that
    $\lambda\cdot\neg$ is $1_{\top}$ follows from the fact that, by
    the interchange law, 
    \begin{align*}
      \lambda\cdot\neg & = (\varphi*\lambda)\cdot\zeta.
    \end{align*}
    Moreover, $(\varphi*\lambda)=(\varphi*\top)$ since
    $(i*\lambda)=1_{U}$.  Thus, $(\lambda\cdot\neg)=1_{\top}$ by another
    application of the interchange law and the fact that
    $(\varphi\cdot\varphi^{-1})=1_{I}$.
  \end{proof}
\end{lemma}
\begin{remark}
  We mention a further equivalence, which we will not require and
  which is easily verified using representability of $I$.
  Namely, $I$ is invertible if and only if it is a Boolean algebra.
\end{remark}
By Lemma \ref{lemma:i_equiv} it follows that, for example,
$\mathbf{2}$ is not a Hopf interval in $\cat$.  Nonetheless, when
$\EE$ it is possible to
construct in the expected manner the free Hopf interval $J$ generated
by $I$.  Specifically, where $\mathbf{I}$ is the category described in
Example (3) from Section \ref{section:examples} and where
$\mathbf{I}\cdot-$ denotes the tensor with $\mathbf{I}$ and
exists since $I$ is representable, we set $J\colon =\mathbf{I}\cdot
U$.  It is then immediate that $J$ is
an invertible interval with symmetry map $\sigma_{J}\colon J\to<150>J$
and there exists a morphism of intervals $\iota\colon I\to<150>J$ with
which exhibits $J$ as the free invertible interval generated by $I$ in
the sense that, for any morphism of intervals $\xi\colon I\to<150>H$
with $H$ invertible, there exists a canonical map of intervals
$\bar{\xi}\colon J\to<150>H$ extending $\xi$ along the inclusion
$\iota$.  Additionally, $J$ classifies the invertible 2-cells in the
2-category structure induced by $I$ as described in the following
lemma (the proof of which is routine).
\begin{lemma}\label{lemma:J}
  Given morphisms $f$ and $g$ in $\EE(B,A)$ together with a 2-cell
  $\alpha\colon f\twocell<150>g$ in the 2-category structure induced by $I$,
  $\alpha$ is invertible if and only if there exists a canonical
  extension $\bar{\alpha}\colon B\tensor J\to<150>A$ such that the following
  diagram commutes:
  \begin{align*}
    \bfig
    \Vtriangle/..>`<-`<-/<300,300>[B\tensor J`A`B\tensor I;\bar{\alpha}`B\tensor\iota`\alpha]
    \efig
  \end{align*}
\end{lemma}
Using Lemma \ref{lemma:J} it is possible to construct meet and join
operations on $J$ as such canonical extensions.  For example, using
representability of $I$, a straightforward calculation shows that,
regarded as a 2-cell $\bottom_{J}\circ i\twocell<150>\iota$, the map
$\iota\circ\wedge\colon I\tensor I\to<150>J$ has as its inverse the vertical
composite $f\cdot(\iota\circ\vee)$ where $f$ is defined to be
\begin{align*}
  I\tensor I\to<350>^{i\tensor I}U\tensor I\to<350>^{\lambda}I\to<350>^{\iota}J\to<350>^{\sigma_{J}}J.
\end{align*}
Thus, by Lemma \ref{lemma:J} there exists a canonical extension
$\barwedge\colon I\tensor J\to<150>J$.  Applying the same trick one more time, using the
symmetry map $\tau\colon I\tensor J\to<150>J\tensor I$, yields the required meet operation
$\wedge_{J}\colon J\tensor J\to<150>J$.  The construction of join is dual.
This construction gives us the following lemma.
\begin{lemma}\label{lemma:J_rep}
  $J$ is a representable interval if $I$ is.
  \begin{proof}
    We have seen that $J$ possesses meet and join operations which
    satisfy the required equations by construction.  Thus, by Theorem
    \ref{theorem:rep_2cat} it suffices to show that if
    $\varphi,\psi\colon B\tensor J\tensor J\two<150>A$ are cells in
    $\EE_{J}(B,A)^{\flat}$ which agree on their boundaries, then they are in fact
    equal.  Since $\iota$ is a morphism of intervals and $I$ is
    representable it follows that $\varphi$ and $\psi$ are equal upon
    precomposition with $(B\tensor\iota)\tensor\iota$.  By
    construction of $J$ it then follows that they are likewise equal
    upon precomposition with $(B\tensor J)\tensor\iota$.  Finally, by
    Lemma \ref{lemma:J}, it follows that $\varphi$ and $\psi$ are
    equal.
  \end{proof}
\end{lemma}
It follows from Lemma \ref{lemma:J_rep} that there also exists an
isofibration model structure on $\EE$ defined with respect to $J$.
However, since only invertible 2-cells feature in the
specification of the isofibration model structure, it follows, by
Lemma \ref{lemma:J}, that this ``$J$-model structure'' coincides with
the original ``$I$-model structure''.  As such, we continue to simply
refer to \emph{the} isofibration model structure on $\EE$ without
reference to either interval $I$ or $J$.  Note though that the
2-category structures induced by $I$ and $J$ \emph{do} in general
differ.  Namely, the 2-cells for $J$ correspond exactly to the invertible
2-cells for $I$.
\begin{proposition}\label{prop:J_hopf}
  $J=(J,\Delta_{J},i_{J},\wedge_{J},\top_{J})$ is the free commutative
  Hopf interval generated by $I$ in the isofibration model structure.
  \begin{proof}
    Since $J$ is an invertible interval it follows from Lemma \ref{lemma:i_equiv} that
    $i_{J}\colon J\to<150>U$ is a weak equivalence.  To see that
    $[\bottom_{J},\top_{J}]\colon U+U\to<150>J$ is a cofibration we first observe
    that, because $U$ is cofibrant , rephrasing the usual argument in terms of the present
    setting, any weak equivalence $p\colon E\to<150>B$ is ``full'' in the
    sense of possessing the right-lifting property with respect to the
    map $[\bottom,\top]\colon U+U\to<150>I$.  Moreover, combined with the fact
    that any $h\colon I\to<150>E$ is invertible (as a 2-cell) whenever
    $p\circ h\colon I\to<150>B$ is, gives by Lemma \ref{lemma:J} that any
    weak equivalence $p$ has the right-lifting property with respect
    to $[\bottom_{J},\top_{J}]$ as well.  

    Observe that $(U+U)$ is the initial object in $\intervalsin(\EE)$
    and $[\bottom_{J},\top_{J}]$ is the coobject part of the induced
    canonical map $(U+U)\to<150>J$ of intervals.  It therefore
    follows, by the remark following Corollary \ref{cor:hopf}, that
    this is a morphism of Hopf objects.  Similarly,
    $i_{J}$ is a morphism of Hopf objects since it is the coobject
    part of the induced map into the terminal object $U$ in
    $\intervalsin(\EE)$.
    
    Finally, for freeness, suppose given another Hopf interval $H$:
    \begin{align*}
      U+U\to<350>^{[a,b]}H\to<350>^{g}U
    \end{align*}
    together with $\xi\colon I\to<150>H$ a morphism of
    commutative Hopf objects. Because $g$ is a weak equivalence,
    there exists an arrow $g'\colon U\to<150>H$
    together with an invertible 2-cell $\varphi\colon g'\circ
    g\twocell<150>1_{H}$.  Then the vertical composite
    $(\varphi*a)\cdot(\varphi^{-1}*b)$ in $\EE(U,H)$ is a 2-cell
    $b\twocell<150>a$ which is, by the fact that $g\xi=i$, the inverse
    of $\xi$.  Thus, by Lemma \ref{lemma:J} there exists a canonical
    extension $\bar{\xi}\colon J\to<150>H$ of $\xi$ which is a
    morphism of Hopf objects and commutes with $[a,b]$ and $g$.
    Finally, $\bar{\xi}$ is trivially the canonical map with these
    properties.
  \end{proof}
\end{proposition}
Berger and Moerdijk \cite{Berger:AHTO} have shown that the existence
of a commutative Hopf interval is one of several conditions which
allow one to lift a model structure from a symmetric monoidal closed
category $\EE$ to the category of reduced operads over $\EE$.
\begin{theorem}[Berger and Moerdijk]
  \label{theorem:Berger_Moerdijk}
  If a symmetric monoidal closed category $\EE$ is a monoidal model
  category such that $\EE$ is cofibrantly generated with cofibrant
  tensor unit $U$, $\EE/U$ has a symmetric monoidal fibrant
  replacement functor and $\EE$ possesses a commutative Hopf interval,
  then there exists a cofibrantly generated model structure on the
  category of reduced operads in which the weak equivalences and
  fibrations are ``pointwise''.
\end{theorem}
Here that the model category is \emph{monoidal} means
that the appropriate internal form of the condition from Definition
\ref{def:model_cat_cat} is satisfied (cf.~\cite{Hovey:MC}).  As such,
in light of Theorem \ref{theorem:Berger_Moerdijk}, we obtain the
following corollary to Proposition \ref{prop:J_hopf}:
\begin{corollary}
  \label{cor:operads}
  Assume $\EE$ possesses a representable interval
  $I$, then when the isofibration model structure is cofibrantly
  generated there is a model structure on the category of reduced
  operads over $\EE$ in which the fibrations and weak equivalences are
  pointwise.
  \begin{proof}
    By the fact that all objects in the isofibration model structure
    are fibrant and Proposition \ref{prop:J_hopf} it remains, in order
    to be able to apply Theorem \ref{theorem:Berger_Moerdijk}, only to
    verify that the model structure is monoidal.  For this, it
    suffices, by the definition of fibrations in the isofibration
    model structure, to note that if $f\colon E\to<150>B$ is a 
    (trivial) fibration, then so is $f_{*}\colon [X,E]\to<150>[X,B]$ for any object
    $X$.  First, that $f_{*}$ is a fibration when $f$ is follows from the
    tensor-hom adjunction.  Next, that $f_{*}$ is
    an equivalence when $f$ is follows from the fact that the map
    \begin{align*}
      \EE(E,B)\to<350>\EE\bigl([X,E],[X,B]\bigr)
    \end{align*}
    which sends an arrow $f\colon E\to<150>B$ to $f_{*}$ is a functor and
    therefore preserves isomorphic 2-cells.
  \end{proof}
\end{corollary}

\subsection*{Acknowledgements}

The main results of this paper occur in the cartesian case in my Ph.D. thesis
\cite{Warren:PhD} --- where I used invertible intervals to obtain models of
intensional Martin-L\"{o}f type theory --- and I
would first and foremost like to thank my thesis supervisor Steve Awodey for
his numerous valuable discussions and suggestions, and for his
comments on a draft of this paper.  This paper has also benefitted
from useful discussions with and comments by Nicola Gambino, Pieter
Hofstra, Tom Leinster, Peter LeFanu Lumsdaine, Alex Simpson, Alexandru Stanculescu and
Thomas Streicher.  I thank the members of the Logic and Foundations of
Computing group at the University of Ottawa for their support during
the preparation of this paper, and for giving me the opportunity to
give a series of talks on this topic.  I am also especially grateful
to the referee for noticing an error in an earlier version of this
paper and for making a number of valuable suggestions including the
suggestion to use the work of Brown and Mosa on double categories with
connections.  Finally, I thank the National Science Foundation for its
support while this paper was being revised.\footnote{This material is
  based upon work supported by the National Science Foundation under
  agreement No. DMS-0635607.  Any opinions, findings and conclusions or
  recommendations expressed in this material are those of the author and
  do not necessarily reflect the views of the National Science Foundation.}

\appendix

\section{Hopf objects and lattices in a symmetric monoidal category}\label{section:monoids_comonoids_hopf}

In this appendix we provide the full definitions of comonoids, Hopf
objects and distributive lattices in a symmetric monoidal category.

\subsection{Monoids, comonoids and Hopf objects}

A \myemph{monoid} $(M,\eta,m)$ in a symmetric monoidal category $\EE$
is given by an object $M$ of $\EE$ together with arrows
$\eta\colon U\to<150>M$ and $m\colon M\tensor M\to<150>M$ satisfying the following
diagrams commute:
\begin{align*}
  \bfig
  \square(0,0)/<-`->`->`->/<600,400>[M\tensor U`M`M\tensor
  M`M;\rho^{-1}`M\tensor\eta`M`m]
  \square(600,0)/->`->`->`<-/<600,400>[M`U\tensor M`M`M\tensor M;\lambda^{-1}`
  `\eta\tensor M`m]
  \efig
\end{align*}
and
\begin{align*}
  \xy
  {\ar^{\alpha}(0,800)*+{M\tensor(M\tensor
      M)};(1200,800)*+{(M\tensor M)\tensor M}};
  {\ar_{M\tensor m}(0,800)*+{M\tensor(M\tensor
      M)};(0,400)*+{M\tensor M}};
  {\ar^{m\tensor M}(1200,800)*+{(M\tensor M)\tensor
      M};(1200,400)*+{M\tensor M}};
  {\ar_{m}(0,400)*+{M\tensor M};(600,0)*+{M}};
  {\ar^{m}(1200,400)*+{M\tensor M};(600,0)*+{M}};
  \endxy
\end{align*}

A \myemph{comonoid} $(G,\epsilon,\Delta)$ in a symmetric monoidal
category $\EE$ is given by an object $G$ of $\EE$ together with arrows
$\epsilon\colon G\to<150>U$ and $\Delta\colon G\to<150>G\tensor G$ such that the
following diagrams commute:
\begin{align*}
  \bfig
  \square(0,0)/<-`->`->`->/<600,400>[M\tensor M`M`U\tensor
  M`M;\Delta`\epsilon\tensor M`M`\lambda]
  \square(600,0)/->`->`->`<-/<600,400>[M`M\tensor M`M`M\tensor U;\Delta`
  `M\tensor\epsilon`\rho]
  \efig
\end{align*}
and
\begin{align*}
  \xy
  {\ar_{\Delta}(600,800)*+{M};(0,400)*+{M\tensor M}};
  {\ar^{\Delta}(600,800)*+{M};(1200,400)*+{M\tensor M}};
  {\ar_{\Delta\tensor M}(0,400)*+{M\tensor M};(0,0)*+{(M\tensor
      M)\tensor M}};
  {\ar^{M\tensor\Delta}(1200,400)*+{M\tensor
      M};(1200,0)*+{M\tensor(M\tensor M)}};
  {\ar_{\alpha^{-1}}(0,0)*+{(M\tensor M)\tensor
      M};(1200,0)*+{M\tensor(M\tensor M)}};
  \endxy
\end{align*}

A \myemph{(commutative) Hopf object} in $\EE$ is a structure
$(H,\eta,m,\epsilon,\Delta)$ such that $(H,\epsilon,\Delta)$ is a
comonoid, $(H,\eta,m)$ is a (commutative) monoid, and the maps $m$ and $\eta$ are
comonoid homomorphisms (cf.~\cite{Berger:AHTO}).  Here note that
$H\tensor H$ is given the structure of a comonoid via the map,
constructed using the symmetry $\tau$, which (schematically) sends
$x\tensor y$ to $\bigl((x\tensor y)\tensor(x\tensor y)\bigr)$.  

\subsection{Lattices}

Assume that $(L,\epsilon,\Delta)$ is a comonoid in $\EE$.  Then $L$ is
a \myemph{lattice} if there are maps $\vee\colon L\tensor L\to<150>L$ and
$\wedge\colon L\tensor L\to<150>L$ such that both $\vee$ and $\wedge$ are
associative, commutative, and following diagrams commute:
\begin{align}\label{eq:idempotent}
  \bfig
  \square<600,400>[L`L\tensor L`L\tensor
  L`L;\Delta`\Delta`\vee`\wedge]
  \morphism(0,400)<600,-400>[L`L;L]
  \efig
\end{align}
and
\begin{align}\label{eq:absorption}
  \xy
  {\ar^{\Delta\tensor L}(0,700)*+{L\tensor L};(800,700)*+{(L\tensor
      L)\tensor L}};
  {\ar^{\alpha^{-1}}(800,700)*+{(L\tensor L)\tensor
      L};(1600,700)*+{L\tensor(L\tensor L)}};
  {\ar_{L\tensor\epsilon}(0,700)*+{L\tensor L};(0,300)*+{L\tensor U}};
  {\ar_{\rho}(0,300)*+{L\tensor U};(800,0)*+{L}};
  {\ar^{L\tensor\clubsuit}(1600,700)*+{L\tensor(L\tensor
      L)};(1600,300)*+{L\tensor L}};
  {\ar^{\diamondsuit}(1600,300)*+{L\tensor L};(800,0)*+{L}};
  \endxy
\end{align}
for $\clubsuit=\vee$ and $\diamondsuit=\wedge$, or $\clubsuit=\wedge$
and $\diamondsuit=\vee$.  

A lattice $L$ is \myemph{unital} if there exist maps
$\top,\bottom\colon U\two<150>L$ such that $\top$ is a unit for $\wedge$,
$\bottom$ is a unit for $\vee$, and the following diagram commutes:
\begin{align}\label{eq:unital_lattice}
  \bfig
  \square(0,0)<600,400>[U\tensor L`L\tensor L`L`L;t\tensor
  L`\lambda`\diamondsuit`t\circ\epsilon]
  \square(600,0)/<-`->`->`<-/<600,400>[L\tensor L`L\tensor
  U`L`L;L\tensor t``\rho`t\circ\epsilon]
  \efig
\end{align}
for $\diamondsuit=\wedge$ and $t=\bottom$, or $\diamondsuit=\vee$ and $t=\top$.

A lattice $L$ is distributive if the further diagram commutes:
\begin{align}\label{eq:dist_lattice}
  \xy
  {\ar^{L\tensor\clubsuit}(600,1000)*+{L\tensor(L\tensor
      L)};(1600,1000)*+{L\tensor L}};
  {\ar^{\diamondsuit}(1600,1000)*+{L\tensor L};(2200,700)*+{L}};
  {\ar_{\alpha\circ(\Delta\tensor (L\tensor
      L))}(600,1000)*+{L\tensor(L\tensor L)};(0,700)*+{\bigl((L\tensor
      L)\tensor L\bigr)\tensor L}};
  {\ar_{\alpha^{-1}\tensor L}(0,700)*+{\bigl((L\tensor L)\tensor
      L\bigr)\tensor L};(0,300)*+{\bigl(L\tensor(L\tensor
      L)\bigr)\tensor L}};
  {\ar_{\tau\tensor L}(0,300)*+{\bigl(L\tensor(L\tensor
      L)\bigr)\tensor L};(600,0)*+{\bigl((L\tensor L)\tensor
      L\bigr)\tensor L}};
  {\ar_{\alpha^{-1}}(600,0)*+{\bigl((L\tensor L)\tensor
      L\bigr)\tensor L};(1600,0)*+{(L\tensor L)\tensor
      (L\tensor L)}};
  {\ar_{\diamondsuit\tensor\diamondsuit}(1600,0)*+{(L\tensor L)\tensor
      (L\tensor L)};(2200,300)*+{L\tensor L}};
  {\ar_{\clubsuit}(2200,300)*+{L\tensor L};(2200,700)*+{L}};
  \endxy
\end{align}
for either  $\clubsuit=\vee$ and $\diamondsuit=\wedge$, or
$\clubsuit=\wedge$ and $\diamondsuit=\vee$.
%% 
%% BIBLIOGRAPHY INFO
%%
\nocite{CRM:ACSMHC}
\bibliographystyle{amsplain}
\bibliography{intervals}
%%
%% END DOCUMENT
%%
\end{document}